\documentclass[11pt]{article}
\usepackage[T1]{fontenc}
\usepackage{lmodern,microtype}
\usepackage[margin=1in]{geometry}
\usepackage{amsmath,amsxtra,amssymb,amscd,amsthm,mathtools}
\usepackage{bbm,mathrsfs,xparse}
\usepackage{physics}
\usepackage[backend=biber,style=alphabetic,sorting=ynt,giveninits=true,isbn=false]{biblatex}
\DeclareFieldFormat{url}{%
  \iffieldundef{usera}{\url{#1}}{\href{#1}{\printfield{usera}}}}
\usepackage{booktabs,needspace}
\usepackage[dvipsnames]{xcolor}
\usepackage{tikz}
\usetikzlibrary{arrows.meta,positioning,calc,decorations.pathreplacing}
\definecolor{figNavy}{RGB}{35,58,80}
\definecolor{figTeal}{RGB}{63,116,110}
\definecolor{figOchre}{RGB}{164,120,53}
\definecolor{figGray}{RGB}{126,135,141}
\definecolor{figLight}{RGB}{224,229,232}
\definecolor{figWash}{RGB}{248,241,226}
\colorlet{darkblue}{MidnightBlue}
\colorlet{darkgreen}{ForestGreen}

\newtheorem{thm}{Theorem}[section]

\newtheorem{cor}[thm]{Corollary}
\newtheorem{lem}[thm]{Lemma}

\newtheorem{prop}[thm]{Proposition}

\theoremstyle{remark}
\newtheorem{remark}[thm]{Remark}
\theoremstyle{definition}

\newcommand{\bean}{\begin{eqnarray}}
\newcommand{\eean}{\end{eqnarray}}
\newcommand{\be}{\begin{displaymath}}
\newcommand{\ee}{\end{displaymath}}
\newcommand{\bea}{\begin{eqnarray*}}
\newcommand{\eea}{\end{eqnarray*}}

\newcommand{\thmref}[1]{Theorem~\ref{#1}}
\newcommand{\secref}[1]{Section~\ref{#1}}
\newcommand{\lemref}[1]{Lemma~\ref{#1}}
\newcommand{\propref}[1]{Proposition~\ref{#1}}
\newcommand{\corref}[1]{Corollary~\ref{#1}}

\newcommand{\nc}{\newcommand}

\makeatletter
\newcommand{\oset}[3][0ex]{%
  \mathrel{\mathop{#3}\limits^{
    \vbox to#1{\kern-2\ex@
    \hbox{$\scriptstyle#2$}\vss}}}}
\makeatother
\newcommand{\brac}[1]{\left( #1 \right)}                                % brackets
\newcommand{\sbrac}[1]{\left\{ #1 \right\}}                             % set
\newcommand{\abrac}[1]{\left\langle #1\right\rangle}                    % angled

\newcommand{\Exponential}{\operatorname{Exp}}

\newcommand{\E}[1]{\mathbb{E}\left[#1\right]} 							% expectation
\newcommand{\sE}[2]{\mathbb{E}_{#2}\left[#1\right]}					    % expectation w/ sub

\newcommand{\Var}[1]{\mathbf{Var}\left(#1\right)}					    % variance

\newcommand{\Prob}[1]{\mathbb{P}\left(#1\right)}						% probability
\newcommand{\I}[1]{{\mathbf 1}_{#1}}									% indicator
\newcommand{\Esub}[2]{{\mathbb E_{#1}}\left[#2\right]}

\newcommand{\iid}{\text{i.i.d.}}
\newcommand{\as}{\text{a.s.}}
\DeclareDocumentCommand{\Pwto} {o} {
	\IfNoValueTF {#1}
	{\overset{n\to\infty}{\longrightarrow}}
	{ \xrightarrow[ #1 \to \infty]{\Pr }}
}
\DeclareDocumentCommand{\Prto} {o} {
	\IfNoValueTF {#1}
	{\overset{\Pr}{\longrightarrow}}
	{ \xrightarrow[ #1 \to \infty]{\Pr }}
}
\DeclareDocumentCommand{\Asto} {o} {
	\IfNoValueTF {#1}
	{\overset{\operatorname{a.s.}}{\longrightarrow}}
	{
		\xrightarrow[ #1 \to \infty]{\operatorname{a.s.} }
	}
}
\DeclareDocumentCommand{\Mgfto} {o} {
	\IfNoValueTF {#1}
	{\overset{\operatorname{mgf}}{\longrightarrow}}
	{ \xrightarrow[ #1 \to \infty]{\operatorname{mgf} }}
}

\DeclareDocumentCommand{\Wkto} {o} {
	\IfNoValueTF {#1}
	{\overset{(d)}{\longrightarrow}}
	{ \xrightarrow[ #1 \to \infty]{(d) }}
}

\DeclareDocumentCommand \LPto { O{1} }
{\overset{\operatorname{L^{#1}}}{\longrightarrow}}
\newcommand\restr[2]{{% we make the whole thing an ordinary symbol
  \left.\kern-\nulldelimiterspace % automatically resize the bar with \right
  #1 % the function
  \vphantom{\big|} % pretend it's a little taller at normal size
  \right|_{#2} % this is the delimiter
  }}

\nc{\on}{\operatorname}
\nc{\ch}{\mbox{ch}}
\nc{\Z}{{\mathbb Z}}
\nc{\C}{{\mathbb C}}
\nc{\prob}{\mathbb{P}}

\nc{\indic}{\mathbbm{1}}
\nc{\indistr}{\overset{d}{\longrightarrow}}
\nc{\inprob}{\overset{p}{\longrightarrow}}
\nc{\N}{{\mathbb N}}
\nc{\pone}{{\mathbb C}{\mathbb P}^1}
\nc{\pa}{\partial}
\nc{\F}{{\mathcal F}}
\nc{\arr}{\rightarrow}
\nc{\larr}{\longrightarrow}
\nc{\al}{\alpha}
\nc{\ri}{\rangle}
\nc{\lef}{\langle}
\nc{\W}{{\mathbb W}}
\nc{\la}{\lambda}
\nc{\lapl}{\mathcal{L}}
\nc{\ep}{\lambda}
\nc{\su}{\widehat{{\mathfrak s}{\mathfrak l}}_2}
\nc{\sw}{{\mathfrak s}{\mathfrak l}}
\nc{\g}{{\mathfrak g}}
\nc{\h}{{\mathfrak h}}
\nc{\n}{{\mathfrak n}}
\nc{\G}{\widehat{\g}}
\nc{\De}{\Delta}
\nc{\gt}{\widetilde{\g}}
\nc{\Ga}{\Gamma}
\nc{\one}{{\mathbf 1}}
\nc{\z}{{\mathfrak Z}}
\nc{\La}{\Lambda}
\nc{\wt}{\widetilde}
\nc{\wh}{\widehat}
\nc{\cri}{_{\kappa_c}}
\nc{\si}{\sigma}
\nc{\el}{\ell}
\nc{\bi}{\bibitem}
\nc{\om}{\omega}
\nc{\ol}{\overline}
\nc{\dzz}{\frac{dz}{z}}
\nc{\mc}{\mathcal}
\nc{\Cal}{\mathcal}
\nc{\bb}{{\mathfrak b}}
\nc{\ot}{\otimes}
\nc{\R}{{\mathbb R}}
\nc{\yy}{{\mc Y}}
\nc{\ga}{\gamma}
\nc{\us}{\underset}
\nc{\opl}{\oplus}
\nc{\Fq}{{\mathbb F}_q}
\nc{\Mq}{{\mathcal M}}
\nc{\Rep}{\on{Rep}}
\nc{\sssec}{\subsubsection}
\nc{\ssec}{\subsection}
\nc{\lan}{\langle}
\nc{\ran}{\rangle}

\newcommand\cA{\mathcal A}
\newcommand\cB{\mathcal B}
\newcommand\cC{\mathcal C}
\newcommand\cD{\mathcal D}
\newcommand\cE{\mathcal E}
\newcommand\cF{\mathcal F}

\newcommand\cH{\mathcal H}

\newcommand\cL{{\mathcal L}}

\newcommand\cN{\mathcal N}

\newcommand\cT{{\mathcal T}}

\newcommand\cX{{\mathcal X}}

\newcommand\cZ{{\mathcal Z}}

\nc{\D}{\mathcal D}
\nc{\Vect}{\on{Vect}}
\nc{\ghat}{\G}
\nc{\T}{\mc T}
\nc{\Tloc}{\T^\g_{\on{loc}}}
\nc{\vac}{|0\ran}
\nc{\Wick}{{\mb :}}
\nc{\mb}{\mathbf}
\nc{\delz}{\partial_z}
\nc{\K}{{\cali K}}
\nc{\cali}{\mathcal}
\nc{\li}{\mathfrak l}
\nc{\lt}{\widetilde{\li}}
\nc{\astar}{a^*}
\nc{\ka}{\kappa}

\nc{\OO}{{\mc O}}
\nc{\AutO}{\on{Aut}\OO}
\nc{\DerO}{\on{Der}\OO}
\nc{\DerpO}{\on{Der}_+\OO}
\nc{\Au}{{\mc A}ut}
\nc{\mf}{\mathfrak}
\nc{\V}{{\mc V}}
\nc{\hh}{\wh{\h}}

\nc{\pp}{{\mathfrak p}}
\nc{\mm}{{\mathfrak m}}
\nc{\rr}{{\mathfrak r}}
\nc{\gr}{\on{gr}}
\nc{\Spe}{\on{Spec}}
\nc{\rv}{\rho^\vee}
\nc{\can}{\on{can}}
\nc{\CC}{\on{Op}_G(D))}
\nc{\Op}{\on{Op}_G(D)}
\nc{\MOp}{\on{MOp}_G(D)}
\nc{\Db}{{\mathbb D}}
\nc{\ww}{w}

\nc{\oQl}{\ol{{\mathbb Q}}_\ell}
\nc{\oFq}{\ol{{\mathbb F}}_q}
\nc{\Q}{{\mathbb Q}}
\nc{\Ql}{{\mathbb Q}_\ell}
\nc{\bs}{\backslash}
\nc{\AD}{{\mathbb A}}
\nc{\M}{{\mc M}}
\nc{\Bun}{\on{Bun}}
\nc{\hl}{h^{\leftarrow}}
\nc{\hr}{h^{\rightarrow}}
\nc{\supp}{\on{supp}}
\nc{\He}{\on{H}}
\nc{\Aut}{\on{Aut}}
\nc{\Ll}{{\mc L}}
\nc{\Coh}{{{\mathcal C}oh}}
\nc{\ovc}{\overset{\circ}}
\nc{\Hav}{\on{H}}
\nc{\Mod}{\on{Mod}}
\nc{\kk}{{\mathfrak k}}
\nc{\vf}{\varphi}
\nc{\Gr}{\on{Gr}}
\nc{\gen}{\on{gen}}
\nc{\IC}{\on{IC}}
\nc{\Jac}{\on{Jac}}

\usepackage[colorlinks=true,linkcolor=blue,citecolor=blue,urlcolor=blue]{hyperref}
\hypersetup{pdftitle={Free energy of bipartite random assignment},pdfstartview=Fit}
\renewcommand{\dd}{\,\mathrm d}
\newcommand{\weight}{\operatorname{wt}}
\makeatletter
\newcommand{\eqlabelalias}[1]{\ltx@label{#1}}
\makeatother
\title{Free energy of the bipartite random assignment problem}
\author{Jacob Shkrob\thanks{Courant Institute of Mathematical Sciences, New York University.
  \href{mailto:jas10184@nyu.edu}{\texttt{jas10184@nyu.edu}}}}
\date{\today}
\begin{document}
\maketitle
\begin{abstract}
We compute the limiting free energy of the bipartite random assignment problem for every pseudodimension $d>0$ and every inverse temperature $\la > 0$ at the scale $n^{1/d}$. We show that the limit is given by a rigorous replica-symmetric formula determined by the unique solution of a recursive distributional equation on a Poisson-weighted infinite tree, confirming conjectures of Aldous and Talagrand. As a consequence, for a fixed inverse temperature $\la > 0$, we deduce exponential concentration of the cost of a Gibbs assignment and determine the exponential growth rate of the number of near-optimal assignments.

% Random assignment pairs~$n$ row vertices with~$n$ column vertices, using each vertex once. We prove a limiting free-energy formula when the independent nonnegative costs satisfy~$\Prob{C\le s}\sim\ka s^d$ near zero, for any~$d>0$. At inverse temperature of order~$n^{1/d}$, a unique recursive law on a Poisson tree determines the limit. The proof first allows unmatched vertices, where deletion gives a bounded local recursion. A weighted completion estimate then restores the perfect-assignment constraint. Logarithmic tail bounds and an endpoint comparison identify the unique limiting message and justify the pressure formula. We give the definitions, comparison estimates, and order of limits needed for this argument.
\end{abstract}

\section{Introduction}\label{sec:introduction}
The assignment problem in the complete bipartite graph $K_{n,n}$ is a well-known problem in combinatorial optimization: given two sets $L_n = \sbrac{1_L,\dots,n_L}$ and $R_n = \sbrac{1_R, \dots, n_R}$ of $n$ vertices, which we call the column and row vertices respectively, and a cost matrix $C = (C_{i,j}) \in \R^{n\times n}$ such that every edge is of the form $(j_L,i_R)$ with cost $C_{i,j}$, find the perfect matching $\pi \in S_n$ or \textit{assignment} which minimizes 
$$
    \min_{\pi \in S_n} \sum_{i=1}^n C_{i,\pi(i)}.
$$
Much has been studied about the limiting value of the optimal assignment under independent and identically distributed costs $C_{i,j}$. Aldous \cite{Aldous} famously resolved the limiting value of the assignment problem for $C_{i,j} \sim \Exponential(1)$, showing that in the limit the expected cost of the minimum assignment converges to $\zeta(2) = {\pi^2}/{6}$. Furthermore, W\"astlund \cite{Wastlund} extended this limit for cost distributions in pseudodimension $d \geq 1$ with CDF given by 
\begin{equation}\label{dim:F}
    \Prob{C \leq s} \sim \kappa s^{d}, \quad s\downarrow 0,
\end{equation}
for some $\kappa > 0$, finding that for the complete graph the minimal assignment cost rescaled by $n^{-(1-1/d)}$ converges to a deterministic constant depending on $d$.\footnote{The term \textit{pseudodimension} comes from the mean-field perspective of the matching problem under random distances: typical distances scale like $s^d$ for small enough $s$.} In particular, one recovers Aldous' exponential cost model for the complete graph $K_n$ from W\"astlund's proof by setting $\kappa=d=1$. Both of these results circumvent the non-rigorous approach from statistical physics which treats the cost $\sum_{i=1}^n C_{i,\pi(i)}$ as a Hamiltonian and study the low temperature $(\la \to \infty)$ limit of the following rescaled Gibbs measure $\mu_{\la,n}$ over the space of perfect matchings $M \in \mathscr{M}_n$:
$$
    \mu_{\la, n}(M) \propto {\exp({- \la (\kappa n)^{1/d} \sum_{e \in M} C_{e}})}.\footnote{Here, we use rescaling $x_{i,j} = (\kappa n)^{1/d}C_{i,j}$ to ensure that the number of low-cost edges for each vertex is of order one: for fixed $x$, 
    $$
        n \Prob{x_{1,1}\leq x} = n F(x/(\kappa n)^{1/d}) \Pwto x^{d}.
    $$}
$$
Throughout, we will refer to the rescaled weights $x_e=x_{i,j}=(\kappa n)^{1/d}C_{i,j}$ for each edge $e=(j_L,i_R)$ and use the function $w_{\la}(x)=e^{-\la x}$ to denote the per-edge contribution to the Hamiltonian. An important quantity in statistical physics is the \textit{free energy} of the associated Gibbs measure $\mu_{\la,n}$ in question. For a fixed matching $M$, we define the weight of $M$ as $\weight_{\la}(M) = \prod_{e \in M} w_{\la}(x_e)$. Writing $\mathscr{M}_n$ for the set of perfect matchings of $K_{n,n}$, the partition function and pressure are defined to be the following quantities
$$
    \cZ_n(\la) := \sum_{M \in \mathscr{M}_n} \weight_\la(M), \qquad p_n(\la) := \frac{1}{n} \log \cZ_n(\la).
$$
Moreover, we define the free energy of the Gibbs measure $\mu_{\la,n}$ as $\E{p_n(\la)}$. As noted by Talagrand \cite{talagrand2003assignment} for the uniform cost case, it was not clear a priori how to extend previously known techniques from the Sherrington-Kirpatrick, Ising perceptron, or $k$-SAT models to study the disordered system associated to the assignment problem. In particular, nothing was known about this model even for a small region of $\la < \la_0$ for some $\la_0 > 0$, except in the case when $|L_n| = n$ and $|R_n| = \lfloor (1 + \alpha)n\rfloor$ (i.e. the number of vertices of the bipartite graph differs asymptotically by a fixed capacity $\alpha$ \cite{talagrand2003assignment}) and progress on this question was subsequently listed as an important open problem by Aldous for many years \cite{Aldous}. Non-rigorous predictions of replica symmetry have been known for all temperatures $\la > 0$ are well-known in the literature since the results of M\'ezard and Parisi \cite{mezard1985replicas, mezard1987solution}. 
\section*{Main Results}

Our main result settles the limiting pressure and overlap formula for the matching problem in every pseudo-dimension $d>0$. In fact, our proof completely bypasses the need for a capacity at all and as a result extends the replica-symmetric phase for all inverse temperatures $\la > 0$. Before we state the main theorem, we need a few definitions. We define $\xi = (\xi_i)_{i\geq 1}$ listed in increasing order ($\xi_1 < \xi_2 < \dots$) to be the Poisson point process on $(0,\infty)$ with intensity measure $\nu_d(\mathrm{d}x)$ given by
\begin{equation*}
 \nu_d(\mathrm{d} x)=d x^{d-1}\dd x,
\end{equation*}
so that $\nu_d((0,x])=x^d\,.$ In particular, the number of points of $\xi$ in a given Borel set~$A\subseteq \R$ with~$\nu_d(A)<\infty$ has a Poisson distribution with mean~$\nu_d(A)$ \cite{LastPenrose}. Throughout, we assume that the cost matrix for $n$ vertices comes from an infinite array $\mathbf{C} = (C_{i,j})_{i\geq 1, j\geq 1}$ of independent copies of the random variable $C$ whose CDF is given by \eqref{dim:F}.

\begin{thm}[Free-energy limit]\label{dim:main}\label{thm:main}
Suppose that the law of each cost satisfies \eqref{dim:F} for some $d>0$ and $\kappa > 0$, and suppose $\xi$ is a Poisson point process with intensity measure $\nu_d(dx)$. Then for every~$\la>0$, there exists a random variable $\cH=\cH_\la$ on~$(0,\infty)$, unique in law, satisfying the following equality in distribution for $\iid$ copies $(\cH_i)_{i\geq 1}$ of $\cH$ independent of $\xi$:
\begin{equation}\label{dim:rde}\eqlabelalias{eq:hard}
 \cH \stackrel{\mathrm d}=\brac{\sum_{i\ge1}e^{-\la\xi_i}\cH_i}^{-1}\,.
\end{equation}
Furthermore, the pressure $p_n(\la) = \frac{1}{n}\log \cZ_n(\la)$ converges almost surely as $n\to \infty$ to the deterministic limit
\begin{equation}\label{dim:formula}\eqlabelalias{eq:pressure}
 p(\la)=-2\E{\log \cH}-
 \int_0^\infty\E{\log(1+e^{-\la x}\cH\cH^\prime)}\,\nu_d(\mathrm{d}x),
\end{equation}
where~$\cH^\prime$ is an independent copy of~$\cH$. In addition, if the cost distribution in \eqref{dim:F} also satisfies $\E{C^q}<\infty$ for some~$q>0$, the sequence of pressures $(p_n(\la))_{n\geq 1}$ is eventually integrable and the convergence above also holds in~$L^1$.
\end{thm}
\begin{remark}
    Equation \eqref{dim:rde} is known as a \emph{recursive distributional equation} specifying a particular probability law and we call such $\cH$ a \textit{cavity message}. Any solution on~$(0,\infty)$ is required to make the infinite sum positive and finite almost surely. In \secref{sec:soft}, we explain why~$\cH$ is called a {cavity message} by showing that it arises as the weak limit of a finite-graph cavity ratio of partition functions before and after vertex deletion.
\end{remark}
The free energy formula is, to the best of author's knowledge, the first to appear in the literature on the matching problem at the scale $\la (\kappa n)^{1/d}$. In fact, Theorem \ref{thm:main} leads to many previously unknown results about matchings in every inverse temperature $\la > 0$. As a first corollary, we prove a new concentration result for the minimum assignment at every temperature $\la > 0$ around a deterministic limit $u(\la)$. As an additional corollary, we obtain a new \textit{geometric} statement about the set of near-optimal solutions of the optimal matching problem in terms of the pressure $p(\la)$ (the proof of Corollary \ref{cor:cost-budget} appears in \secref{sec:corol}). For a perfect matching~$M$, we let
\[
 \cC_n(M): =\frac{1}{n}\sum_{e\in M}x_e,
\]
be the rescaled cost associated to a matching $M$.

\begin{cor}[Quenched Gibbs energy concentration]\label{cor:energy}
Under the assumptions of \thmref{dim:main}, let~$\la>0$ be a
differentiability point of~$p$ and define $$u(\la):=-p'(\la) = \int_{0}^\infty x\E{\frac{e^{-\la x}\cH \cH^\prime}{1+e^{-\la x}\cH \cH^\prime}}\nu_d(\mathrm{d}x).$$ The pressure
$p$ is differentiable at all but countably many~$\la>0$. For each~$\varepsilon>0$, there is a deterministic~$c_\varepsilon>0$
such that for almost every cost array $\mathbf{C} = (C_{i,j})_{i\geq 1, j\geq 1}$ and all sufficiently large~$n$,
\[
 \mu_{\la,n}\bigl(|\cC_n(M)-u(\la)|>\varepsilon\bigr)
 \le 2e^{-c_\varepsilon n}.
\]
Moreover, the Gibbs average $\abrac{\cC_n}_{\la,n}:=\Esub{\mu_{\la,n}}{\cC_n(M)}$ converges to~$u(\la)$ for almost every $\mathbf{C}$.
\end{cor}
\begin{proof}[Proof Sketch of \corref{cor:energy}]
Note that $p_n$ is a normalized log-sum-exp of affine functions and hence convex in~$\la$ for every $n\geq 1$.
By the definition of convexity and applying Theorem \ref{thm:main} to the temperatures $\la_0, \la_1,$ and $t\la_0 + (1-t)\la_1$ for some $t \in (0,1)$, the limit $p$ is convex on $(0,\infty)$ and differentiable except at countably many points. (The full cavity formula for $-p^{\prime}(\la)$ is proved rigorously in \lemref{lem:pressure-derivative}).

Fix a differentiability point~$\la$ and write $u=u(\la)$. For $|t|<\la$,
the Gibbs weights gives
\begin{equation}\label{eq: logsum}
    \frac1n\log\sE{e^{nt(\cC_n(M)-u)}}{\mu_{\la,n}}
 =p_n(\la-t)-p_n(\la)-tu.
\end{equation}
The limiting right-hand side is $o(|t|)$ as $t\to0$. For a given $\varepsilon > 0$, choose $h = h_k$ sufficiently small and $n$ large so that \eqref{eq: logsum} is at most $h\varepsilon/2$ at both $t = h$ and $t = -h$ by Theorem \ref{thm:main}. Therefore by convergence of the pressure at $\la$ and $\la \pm h$, we obtain
$$
    \frac{1}{n}\log \Esub{\mu_{\la,n}}{e^{\pm n h (\cC_n(M) - u)}} \leq \frac{h\varepsilon}{2}.
$$
Applying Chernoff's inequality to both $h$ and $-h$ gives the stated concentration inequality
$$
    \mu_{\la,n}\brac{\abs{\cC_n(M) - u}\geq \varepsilon}\leq 2e^{-nh\varepsilon/2}.
$$
For each $k$, convexity of $p_n$ gives
\[
 \frac{p_n(\la)-p_n(\la-h_k)}{h_k}
 \le p_n'(\la)
 \le \frac{p_n(\la+h_k)-p_n(\la)}{h_k}.
\]
Taking $n\to\infty$ and then $k\to\infty$ gives $p_n'(\la)\to p'(\la)$.
Differentiating the finite partition function gives
$p_n'(\la)=-\abrac{\cC_n}_{\la,n}$, and so $\abrac{\cC_n}_{\la,n}$ converges to $u(\la)$ almost surely.
\end{proof}

\begin{cor}\label{cor:cost-budget}
Under the hypotheses of \thmref{dim:main}, define $$\alpha_d^\star := -\lim_{\la \to \infty} \frac{p(\la)}{\la}$$
to be the ground-state minimum of the assignment problem. Then, almost surely, we have that as $n\to\infty$, 
$$\min_{M \in \mathscr{M}_n} \cC_n(M) \longrightarrow \alpha_d^\star.$$
Moreover, for each $\epsilon>0$, as $n\to \infty$
\[
    \frac1n\log\#\{M\in\mathscr{M}_n: \cC_n(M)\le \alpha_d^\star+\epsilon\}
    \longrightarrow
    \inf_{\la>0}\{p(\la)+\la(\alpha_d^\star+\epsilon)\}, \quad \as
\]
% There exists $\la_\epsilon>0$ with
% $-p'(\la_\epsilon)=\alpha_d^\star+\epsilon$. Let $V_{n,\epsilon}$
% be the uniform law on this set of assignments and let $\cQ(\la)$
% be the Gibbs overlap limit from the companion paper.
% For every $\delta>0$ there is a deterministic $c>0$ such that,
% almost surely for all sufficiently large $n$,
% \[
%     V_{n,\epsilon}^{\otimes 2}\!\left(
%         \left|\frac{|M_1\cap M_2|}{n}-\cQ(\la_\epsilon)\right|>\delta
%     \right)\le e^{-cn}.
% \]
\end{cor}

\begin{remark}\label{rem:ground-state-cost}
For pseudimension $d = 1$ and $\kappa = 1$, the ground-state constant
$\alpha_d^\star$ was already given by the constant $\zeta(2)$ ~\cite{Aldous}. Indeed, by Theorem \ref{thm:main}, if $U_\la = \la^{-1}\log \cH_\la$ and $S_\la = U_\la + U^\prime_\la$ be a sum of independent copies of $U_\la$, 
$$
    -\frac{p(\la)}{\la} = 2\E{U_\la} + \E{\int_0^\infty \la^{-1}\log(1 + e^{\la(S_\la -x)})\,\mathrm{d}x}.
$$
Every subsequential limit of $U_\la$ solves the recursion $U \overset{\mathrm{d}}{=}\min_{i}(\xi_i - U_i)$ where $(\xi_i)_{i\geq 1}$ is a rate-one Poisson point process, and by \cite[Lemma 5]{Aldous}, the law of $U$ is logistic. Therefore, in the limit as $\la \to \infty$, we get
$$
    \alpha_d^\star = 2\E{U} + \frac{1}{2}\E{(U + U^\prime)^2_+} = \frac{1}{4}\E{(U+U^\prime)^2} = \frac{1}{2}\Var{U} = \zeta(2). 
$$
\begin{remark}
    Note that $u(\la)\downarrow \alpha_d^\star$ as $\la \to \infty$. Since $p$ is convex, $u = -p^\prime$ is nonincreasing. Therefore, for any fixed $\la_0$,
    $$
        -\frac{p(\la)}{\la} = - \frac{p(\la_0)}{\la} + \frac{1}{\la}\int_{\la_0}^\la u(t)\,dt,
    $$
    and so taking $\la \to \infty$ gives $\lim_{\la\to \infty}u(\la)$ on the RHS and the definition of $\alpha_d^\star$ gives the equality. 
\end{remark}
\end{remark}
\ssec{Related Work}
Rigorous work on the bipartite random assignment problem and its limit can be originally traced to Aldous, who originally formulated the finite-temperature assignment recursion and noted the difficulty extending to the finite $\la$ regime~\cite[Section~5.3]{Aldous}. Using methods similar to the cavity method from the statistical physics literature, Talagrand's approach is similar to the one used in this paper, but fails to reach the entire $\la > 0$ regime \cite{talagrand2003assignment}. The main approach in this paper was heavily inspired by the local-recursions used by Bordenave, Lelarge, and Salez to study the monomer-dimer model \cite{BLS}. Their results concern mainly bounded degree graphs and their matching numbers, whereas our recursion is adapted more for the minimal assignment and includes the random weights in the Gibbs measure. The rigorous approach to replica symmetry of the matching problem in general pseudodimension and additional estimates needed for the alternating-path construction appears originally in the works of W\"astlund \cite{Wastlund}. Recent work of Dey and Terlov \cite{DeyTerlov} studied the equivalent model for more general Gibbs measures associated to random optimization problems, but were only able to study their limits at constant temperature. 
\ssec{Proof Outline}
We now illustrate the general proof outline for the free energy as follows (see Figure \ref{fig:proof-structure} for a concise summary). The remaining proof of Corollary \ref{cor:cost-budget} and additional technical lemmas appear in Section \ref{sec:corol} and Appendix \ref{app:simplex}.
\begin{itemize}
    \item In \secref{sec:soft}, we begin by introducing a monomer activity
    $z>0$ to allow for unmatched vertices. The vertex-deletion recursion
    and local convergence to the Poisson-weighted infinite tree give the limiting cavity
    message $\cH_z$, from which we compute the equivalent formula for the pressure $p(\la,z)$ associated to the Gibbs measure with monomer activity parameter $z > 0$.

    \item In \secref{sec:completion}, we complete partial matchings along
    short, low-cost augmenting paths using a standard argument of growing alternating forests at free vertices. By bounding the total weight of the
    partial matchings that lead to each perfect matching, we show that
    $p_n(\la)\to\lim_{z\downarrow0}p(\la,z)$ almost surely as $n\to\infty$.

    \item In \secref{sec:hardlaw}, we use tail bounds for the random variable $\log\cH_z$ to
    pass to the cavity equation at $z=0$. We prove uniqueness of the limiting random variable $\cH$ as in \eqref{dim:rde} by coupling
    two solutions and exchanging the endpoints of an inserted Poisson
    edge. The same tail bounds justify taking $z\downarrow0$ in the
    pressure formula.
\end{itemize}
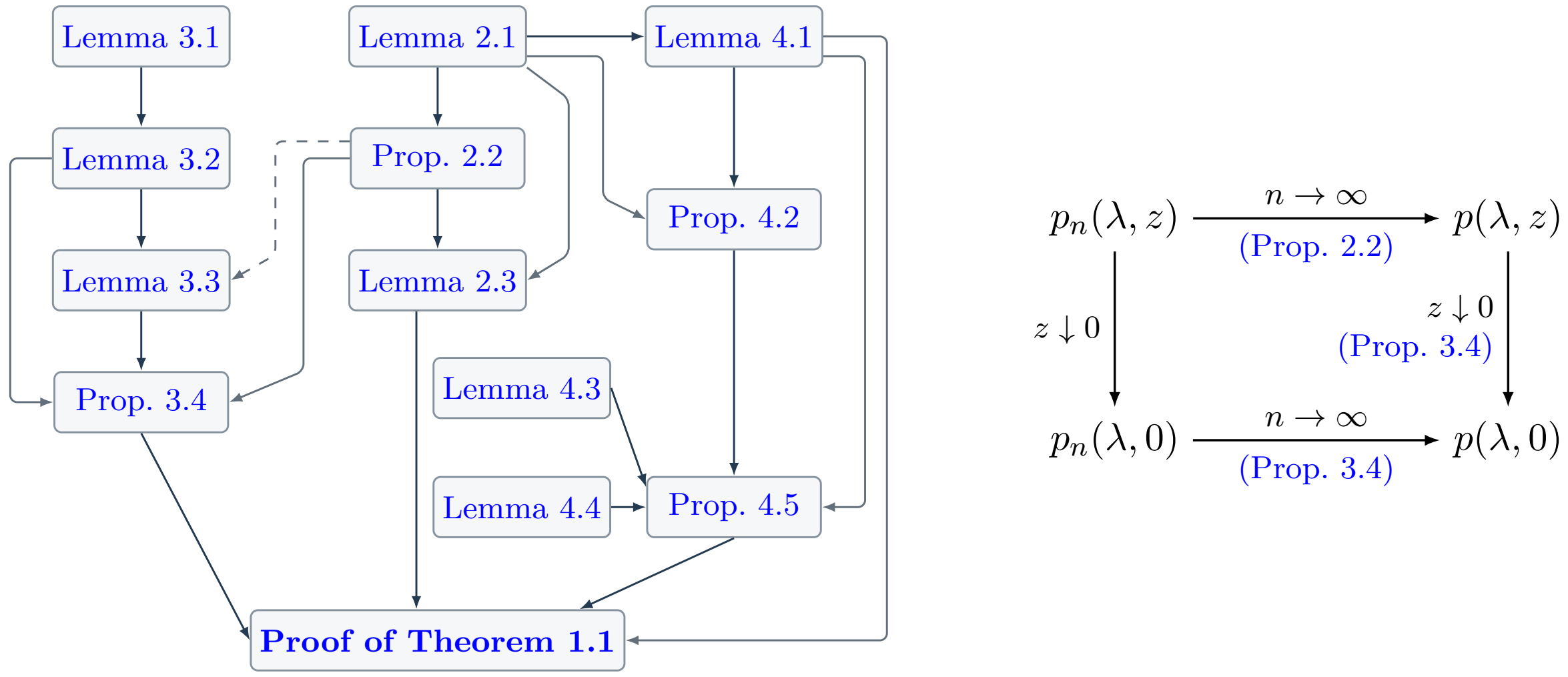
\begin{figure}[!htb]
\centering
\begin{minipage}[c]{.61\linewidth}
\centering
\begin{tikzpicture}[x=.88pt,y=.8pt,
  result/.style={draw=figNavy!55,fill=figNavy!4,rounded corners=2pt,
    line width=.55pt,font=\footnotesize,minimum width=49pt,
    minimum height=17pt,inner xsep=2.5pt,inner ysep=2pt},
  link/.style={-{Latex[length=1.3mm]},draw=figNavy,line width=.55pt,rounded corners=2pt},
  bypass/.style={link,draw=figGray!70!figNavy},
  estimate/.style={bypass,dashed}]
\node[result] (expansion) at (0,0)
  {\hyperref[lem:expander]{Lemma~\ref*{lem:expander}}};
\node[result] (augment) at (0,-43)
  {\hyperref[lem:augment]{Lemma~\ref*{lem:augment}}};
\node[result] (weights) at (0,-86)
  {\hyperref[lem:weights]{Lemma~\ref*{lem:weights}}};
\node[result] (perfect) at (0,-129)
  {\hyperref[prop:perfect]{Prop.~\ref*{prop:perfect}}};

\node[result] (messages) at (95,0)
  {\hyperref[lem:positive]{Lemma~\ref*{lem:positive}}};
\node[result] (positive) at (95,-43)
  {\hyperref[prop:positive]{Prop.~\ref*{prop:positive}}};
\node[result] (bethe) at (95,-86)
  {\hyperref[lem:bethe]{Lemma~\ref*{lem:bethe}}};

\node[result] (tails) at (190,0)
  {\hyperref[dim:log-tails]{Lemma~\ref*{dim:log-tails}}};
\node[result] (existence) at (190,-64.5)
  {\hyperref[prop:existence z = 0]{Prop.~\ref*{prop:existence z = 0}}};
\node[result] (simplex) at (122,-124)
  {\hyperref[lem:simplex]{Lemma~\ref*{lem:simplex}}};
\node[result] (coupling) at (122,-166)
  {\hyperref[lem:invariant-coupling]{Lemma~\ref*{lem:invariant-coupling}}};
\node[result] (unique) at (190,-166)
  {\hyperref[prop:unique]{Prop.~\ref*{prop:unique}}};

\node[result] (formula) at (95,-213)
  {\textbf{\hyperref[sec:assembly]{Proof of Theorem~\ref*{thm:main}}}};

% Completion bounds and the fixed-activity pressure.
\draw[link] (expansion) -- (augment);
\draw[link] (augment) -- (weights);
\draw[link] (weights) -- (perfect);
\draw[bypass] (augment.west) -- (-42,-43) |- (perfect.west);
\draw[estimate] ([yshift=4.8pt]positive.west) -- (43,-37)
  -- (43,-77) -- (weights.east);
\draw[bypass] (positive.west) -- (52,-43)
  -- (52,-118) -- (perfect.east);

\draw[link] (messages) -- (positive);
\draw[link] (positive) -- (bethe);
\draw[bypass] (messages.south east) -- (137,-22)
  -- (137,-77) -- (bethe.east);

% Passage to the zero-activity law.
\draw[link] (messages) -- (tails);
\draw[bypass] ([yshift=-5.6pt]messages.east) -- (148,-7)
  -- (148,-57) -- (existence.west);
\draw[link] (tails) -- (existence);
\draw[link] (existence) -- (unique);
\draw[link] (simplex.east) -- ([yshift=5.6pt]unique.west);
\draw[link] (coupling) -- (unique);
\draw[bypass] ([yshift=-5.6pt]tails.east) -- (232,-7) |- (unique.east);

% Inputs used in the final pressure formula.
\draw[link] (perfect.south) -- (formula.west);
\draw[link] ([xshift=-6pt]bethe.south) -- ([xshift=-6pt]formula.north);
\draw[link] (unique.south) -- ([xshift=40pt]formula.north);
\draw[bypass] (tails.east) -- (239,0) |- (formula.east);
\end{tikzpicture}
\end{minipage}\hfill
\begin{minipage}[c]{.37\linewidth}
\centering
\begin{tikzpicture}[
  link/.style={-{Latex[length=1.5mm]},draw=black,line width=.65pt},
  annotation/.style={font=\footnotesize,text=black}]
\node (softfinite) at (0,0)    {$p_n(\la, z)$};
\node (softlimit)  at (3.9,0)  {$p(\la, z)$};
\node (hardfinite) at (0,-2.2) {$p_n(\la, 0)$};
\node (hardlimit)  at (3.9,-2.2) {$p(\la, 0)$};
\draw[link] (softfinite) --
  node[annotation,midway,above] {$n\to\infty$}
  node[annotation,midway,below]
    {\hyperref[prop:positive]{(Prop.~\ref*{prop:positive})}} (softlimit);
\draw[link] (hardfinite) --
  node[annotation,midway,above] {$n\to\infty$}
  node[annotation,midway,below]
    {\hyperref[prop:perfect]{(Prop.~\ref*{prop:perfect})}} (hardlimit);
\draw[link] (softfinite) --
  node[annotation,midway,left] {$z\downarrow0$} (hardfinite);
\draw[link] (softlimit) --
  node[annotation,midway,left,align=right]
    {$z\downarrow0$\\\hyperref[prop:perfect]{(Prop.~\ref*{prop:perfect})}} (hardlimit);
\end{tikzpicture}
\end{minipage}
\caption{\textbf{Left.} Dependencies in the proof of \thmref{dim:main}.
Arrows point from inputs to results. The dashed arrow marks the row-weight
bound from the proof of \propref{prop:positive} used in \lemref{lem:weights}.
\textbf{Right.} The two orders of limits for the pressure, with
$p_n(\la,0)=p_n(\la)$ and $p(\la,0)=p(\la)$. By \propref{prop:perfect},
the limits agree almost surely, with $L^1$ convergence under the moment condition.}
\label{fig:proof-structure}
\end{figure}

\ssec{Notation}
All logarithms are natural. We write $\stackrel{\mathrm d}=$ for equality in distribution, $\Rightarrow$ for convergence in distribution, $\cL(X)$ for the law of $X$, and $\I{A}$ for the indicator of $A$. As $n\to\infty$, $f_n=O(g_n)$ means $|f_n|\le C|g_n|$ eventually for some constant $C$, and $f_n=o(g_n)$ means $f_n/g_n\to0$. A subscript on $O$ allows the constant to depend on that parameter. We write $f\sim g$ when $f/g\to1$. We use $u^+=u_+=\max(u,0)$ and $u^-=\max(-u,0)$, with $\|X\|_{L^r}=\E{|X|^r}^{1/r}$ and $\Var{X}=\E{(X-\E{X})^2}$. Eventual integrability means~$\E{|p_n(\la)|}<\infty$ for all sufficiently large~$n$; convergence in~$L^1$ means~$\E{|p_n(\la)-p(\la)|}\to0$. The two terms in~\eqref{dim:formula} are finite. The factor~$2$ accounts for the row and column endpoints. The integral corrects for joining two independent endpoint environments by an edge of length~$x$. We say that a family of laws in~$\mathscr{P}((0,\infty))$ is \emph{tight} if for every~$\epsilon>0$ some compact interval~$[a,b]\subset(0,\infty)$ has probability at least~$1-\epsilon$ for every member of the family. A family~$\cX$ of integrable random variables is \emph{uniformly integrable} if
\begin{equation*}
 \lim_{T\to\infty}\sup_{X\in\cX}\E{|X|\I{\{|X|>T\}}}=0\,.
\end{equation*}
For two probability measures $\mu$ and~$\tilde{\mu}$ in $\mathscr{P}((0,\infty))$, a \emph{coupling} is a joint law~$\rho$ on~$(0,\infty)^2$ whose first marginal is $\mu$ and second marginal is $\tilde{\mu}$. We define the set of all couplings as $\Pi(\mu,\tilde{\mu})$. Throughout, if $\cZ_G$ denotes the weighted graph partition function, for an edge $e = (u,v)$, $\cZ_{G-e}$ denotes the partition function deleting only $e$ (keeping $u,v$), $\cZ_{G - \sbrac{u,v}}$ denotes the partition function after deleting $u,v$, together with every edge incident to either vertex, and $\cZ_{G - v}$ denotes the partition function after deleting $v$ and every edge incident to $v$ (including $e$). 

\section{The tree calculation}\label{sec:proof}\label{sec:soft}
We begin with a relaxation of the Gibbs measure by incorporating a monomer activity $z > 0$ to allow for unmatched vertices. This strategy gives a useful recursion for the partition function associated to $z$. We first show that using techniques from local weak convergence, there exists a unique random variable $\cH_z$ that emerges as the limiting cavity message at the root of the PWIT in the infinite-depth limit. We then show that for fixed parameter $z > 0$, the pressure $p_n(\la, z)$ converges almost surely and in $L^1$ to a deterministic limit $p(\la, z)$ depending only on the final limiting message $\cH_z$. Throughout this section, we fix $d,\la>0$ and write $G_n=K_{n,n}$  and $w(x)=e^{-\la x}$.

% Deleting a row vertex leaves unequal class sizes, so the perfect-assignment partition function vanishes on the remaining graph. Allowing unmatched vertices avoids this obstruction: the partition function stays positive after deletion and satisfies an exact recursion. On the limiting tree, deleting the root separates the branches. This is the source of the independent messages in the recursive law. 

\ssec{Vertex deletion and the tree equation}

For a finite bipartite graph~$G$, let~$V(G)$ be its vertex set and let~$\mathscr{M}(G)$ be its set of all possible matchings in $G$, including the empty matching. Write~$G-v$ for the graph obtained by deleting~$v$ and its incident edges, and~$G-\sbrac{u,v}$ for the graph obtained by deleting both vertices $u$ and $v$. Assign to each edge~$e$ the finite nonnegative length~$x_e = (\kappa n)^{1/d}C_e$ and weight~$w_e=w(x_e)$. We now rigorously allow for the introduction of \textit{monomers}, which are isolated vertices in $G$ without any assigned edges, by relaxing the measure $\mu_{\la, n}$ to accommodate \textit{partial matchings}. Let $z > 0$ be the monomer activity parameter (i.e. the additional weight of an unmatched vertex included in the Hamiltonian) and define the new partition function and cavity ratio at $v$ as
\begin{equation*}
\cZ_G(z)=\sum_{M\in\mathscr{M}(G)}z^{|V(G)|-2|M|}\prod_{e\in M}w_e,
 \qquad \cH_{G,z}(v)=\frac{\cZ_{G-v}(z)}{\cZ_G(z)}\,.
\end{equation*}
Similar to the original Gibbs measure $\mu_{\la, n}$, we define $\mu_{G,z}$ to be the Gibbs measure over the space of matchings $\mathscr{M}(G)$ with monomer activity $z > 0$. On $G_n = K_{n,n}$, let $\cZ_{n,k}$ be the sum of matching weights over matchings of size $n-k$, so $\cZ_{n,0}= \cZ_n(\la)$. Therefore, 
\begin{equation}\label{eq:soft-model}
 \cZ_{G_n}(z)=\sum_{k=0}^n z^{2k}\cZ_{n,k}.
\end{equation}
Define the pressure associated to $\cZ_{G_n}(z)$ for $z\ge0$ as
$$
 p_n(\la, z)=n^{-1}\log\cZ_{G_n}(z).
$$
By~\eqref{eq:soft-model}, $p_n(\la, 0)=p_n(\la)$. For $z>0$, the empty matching ensures that every partition function is positive, which implies that $\cH_{G,z}(v)$ is away from $0$. Note that under the law $\mu_{G,z}$, the probability that~$v$ is unmatched is exactly~$z\cH_{G,z}(v)$, and hence $0 < \cH_{G,z}(v) \leq 1/z$ holds for every $v \in V(G)$. For a fixed matching $M \in \mathscr{M}(G)$, $M$ either leaves the vertex $v$ unmatched or pairs it with exactly one of its neighbors $u \sim v$. Therefore, we get the following equality for the partition functions following the famous Heilmann-Lieb recursion relation for monomer-dimer systems \cite{heilmann1972theory}
\begin{equation*}
 \cZ_G(z)=z\cZ_{G-v}(z)
   +\sum_{u\sim v}w_{vu}\cZ_{G-\sbrac{u,v}}(z)\,.
\end{equation*}
Division by~$\cZ_{G-v}(z)$ gives a beautiful recursive formula for the cavity ratio $\cH_{G,z}(v)$ at $v$:
\begin{equation}\label{eq:delete}
 \cH_{G,z}(v)=\brac{z+\sum_{u\sim v}w_{vu}\cH_{G-v,z}(u)}^{-1}\,.
\end{equation}
Throughout, we write $\cH_{n,z}(v) = \cH_{G_n,z}(v)$ dropping the $G_n$-dependence whenever implicit. 

\subsubsection{Poisson-weighted infinite tree and Mecke's equation}
An important object we will use throughout this work is the \emph{Poisson-weighted infinite tree (PWIT) $\cT$}, initially constructed by Aldous \cite{Aldous, sznitman_objective_2004}. To define the infinite tree $\cT$, we first index the vertices by finite words of positive integers, with the empty word~$\varnothing$ at the root. For each vertex~$v \in V(\cT)$, we independently generate a Poisson process~$(\xi_{v,i})_{i\ge1}$ with intensity~$\nu_d$, and connect~$v$ to its child~$vi$ by an edge of weight~$\xi_{v,i}$. Note that there are countably many vertices, and every vertex has infinitely many children. Nevertheless,
\begin{equation*}
 \E{\sum_i w(\xi_{v,i})}=\int_0^\infty e^{-\la x}\nu_d(\dd x)<\infty\,,
\end{equation*}
so the total outgoing weight is therefore finite almost surely for every vertex $v \in V(\cT)$.

\begin{figure}[!htbp]
\centering
\begin{tikzpicture}[x=1cm,y=1cm,font=\small,>=Latex,
 row/.style={rectangle,draw=figNavy,fill=white,line width=.8pt,minimum size=5.3pt,inner sep=0pt},
 col/.style={circle,draw=figTeal,fill=white,line width=.8pt,minimum size=6pt,inner sep=0pt},
 branch/.style={draw=figNavy!65,line width=.65pt}]
 % A schematic thresholded array, with one root row selected.
 \node[anchor=south] at (1.47,3.22) {Scaled cost array};
 \fill[figWash] (0,1.26) rectangle (2.94,1.68);
 \foreach \k in {0,...,7}{
   \draw[figLight,line width=.35pt] (0,.42*\k)--(2.94,.42*\k);
   \draw[figLight,line width=.35pt] (.42*\k,0)--(.42*\k,2.94);
 }
 \draw[figGray!75,line width=.55pt] (0,0) rectangle (2.94,2.94);
 \foreach \x/\y in {.21/2.73,.21/.63,1.47/2.31,1.47/1.05,2.73/.21,1.05/1.89,2.31/.63}
   \fill[figGray] (\x,\y) circle[radius=1.65pt];
 \foreach \x in {.21,1.47,2.73}
   \fill[figOchre] (\x,1.47) circle[radius=2.2pt];
 \node[row] at (-.55,1.47) {};
 \node[anchor=east] at (-.73,1.47) {$j$};
 \draw[figGray,densely dotted] (-.44,1.47)--(-.04,1.47);
 \node[anchor=north,font=\footnotesize] at (1.47,-.18) {dots: $x_{i,j}\le R$};
 % Fixed-depth local limit.
 \draw[-{Latex[length=2mm]},figGray,line width=.65pt] (3.9,1.55)--(5.75,1.55);
 \node[above=4pt] at (4.825,1.55) {$n\to\infty$};
 \node[below=4pt,font=\footnotesize] at (4.825,1.55) {locally};
 \node[anchor=south] at (10.25,3.22) {Poisson neighborhood};
 \coordinate (root) at (10.25,2.9);
 \coordinate (c1) at (7.45,1.55);
 \coordinate (c2) at (10.25,1.55);
 \coordinate (c3) at (13.05,1.55);
 \draw[figOchre,line width=.85pt] (root)--node[midway,above=4pt,inner sep=0pt] {$\xi_1$} (c1);
 \draw[figOchre,line width=.85pt] (root)--node[midway,right=4pt,inner sep=0pt] {$\xi_2$} (c2);
 \draw[figOchre,line width=.85pt] (root)--node[midway,above=4pt,inner sep=0pt] {$\xi_3$} (c3);
 \foreach \p/\x in {c1/6.85,c1/8.05,c2/9.65,c2/10.85,c3/13.05}{
   \draw[branch] (\p)--(\x,.25);
   \draw[branch,densely dotted] (\x,.25)--(\x-.25,-.15);
   \draw[branch,densely dotted] (\x,.25)--(\x+.25,-.15);
   \node[row] at (\x,.25) {};
 }
 \node[row] at (root) {};
 \foreach \p in {c1,c2,c3}{\node[col] at (\p) {};}
 \node[row] at (7.4,-.68) {};
 \node[anchor=west,font=\footnotesize] at (7.58,-.68) {row type};
 \node[col] at (10.3,-.68) {};
 \node[anchor=west,font=\footnotesize] at (10.48,-.68) {column type};
\end{tikzpicture}
\caption{Dots mark lengths~$x_{i,j}\le R$; the shaded row corresponds to the root~$j_R\in R_n$. At fixed depth and cutoff~$R>0$, its neighborhood converges to the Poisson tree with alternating vertex types and independent child processes of intensity~$\nu_d|_{[0,R]}$. Squares denote row type and circles column type. }\label{fig:proof}
\end{figure}
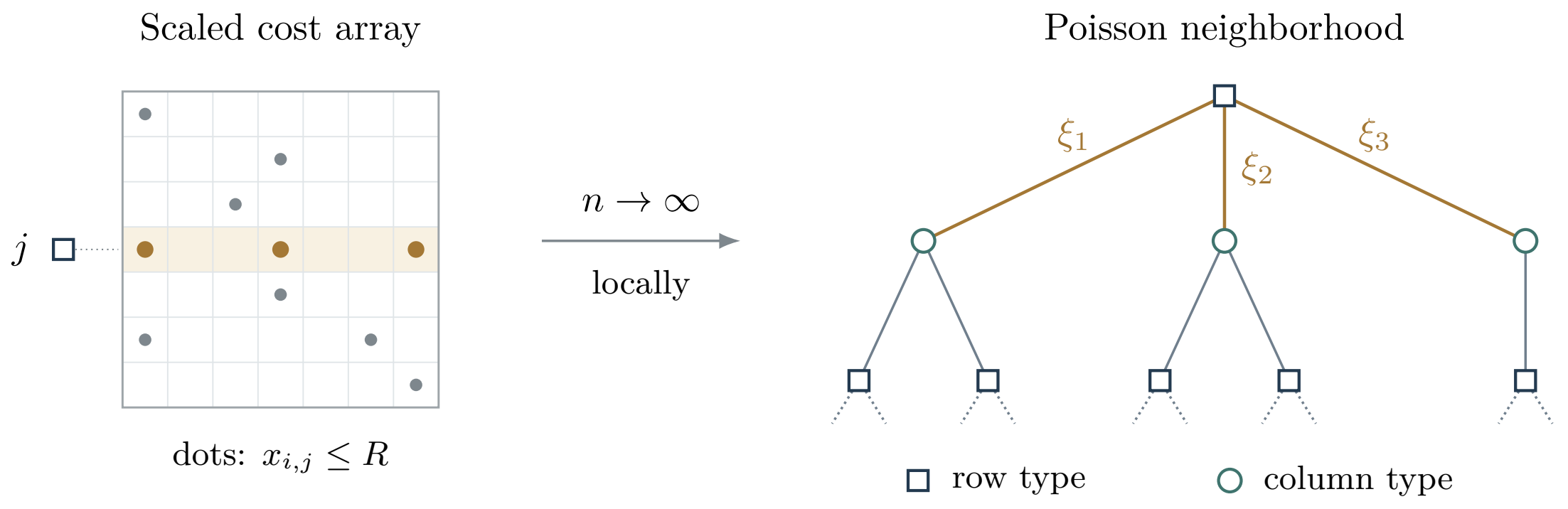

For cutoff $R >0$, define ~$G_{n,R}$ to be subgraph of $G_n$ after deleting edges whose rescaled weights $x_e$ satisfy $x_e > R$, and similarly define $\cT_{R}$
to be the component of the Poisson tree containing the root $\varnothing$ after deleting all
edges with weight larger than $R$. By definition of the Poisson point process with intensity $\nu_d(\mathrm{d}x)$, each vertex of~$\cT_{R}$ has an independent
$\operatorname{Poisson}(R^d)$ number of children, so~$\cT_{R}$ is locally
finite. For fixed~$R$ and~$\el$, an important result is that for a uniformly chosen row or column vertex in $G_{n,R}$, the size-$\el$ neighborhood of the vertex (i.e. the graph induced by all vertices reachable by size-$\el$ paths in $G_{n,R}$) converges in the sense of local weak convergence to the corresponding size-$\el$ neighborhood of the root $\varnothing$ in~$\cT_{R}$. This type of result is the bipartite analogue of W\"astlund's local limit for~$K_n$
\cite[Section~2.2, Proposition~2.2]{Wastlund}. For convenience, we give the exploration argument in the proof of \propref{prop:positive}; see also
Figure~\ref{fig:proof} for an illustration of the truncated graphs $G_{n,R}$ and $\cT_R$.

% To compute the limiting probability that the root is a monomer at fixed
% $z>0$, we first impose a length cutoff~$R$. Local convergence lets us
% evaluate the upper and lower bounds from the depth-$\el$ vertex-deletion
% recursion on a Poisson tree. The two bounds meet as~$\el\to\infty$. 

% To use this local approximation, boundary values far from the root must have a vanishing effect. We will prove this by comparing the two extreme boundary values~$0$ and~$1/z$. First we state the Poisson identity used both here and in the uniqueness argument without monomers.

We also require some standard results from the theory of Poisson point processes, namely the \textit{marking theorem} and \textit{Mecke's equation}. For a Poisson point process $(\xi_i)_{i\geq 1}$ with intensity $\nu_d$, attach independent marks~$Y_i$ that are each distributed according to a probability measure $Q$. By the marking theorem, the counting measure $\cN=\sum_i\delta_{(\xi_i,Y_i)}$, where~$\delta_{(x,y)}$ is the Dirac mass at~$(x,y)$, is a Poisson process with intensity~$\nu_d\otimes Q$~\cite[Theorem~5.6]{LastPenrose}. If~$Y$ is a fresh mark with law~$Q$, independent of~$\cN$, and $\cN - \delta_{(\xi_i, Y_i)}$ is shorthand for the deletion of the point $(\xi_i, Y_i)$ from the point process $\cN$, then for every nonnegative measurable function~$g$,
\begin{equation}\label{eq:mecke}
 \E{\sum_i g\bigl(\xi_i,Y_i,\cN-\delta_{(\xi_i,Y_i)}\bigr)}
 =\int_0^\infty\E{g(x,Y,\cN)}\,\nu_d(\mathrm{d} x)\,.
\end{equation}
Equation \eqref{eq:mecke} is the deletion version of the famous Mecke's equation~\cite[Theorem~4.1]{LastPenrose}. The formula also holds for signed~$g$ once the corresponding absolute-value integral is finite. 

\begin{figure}[!htbp]
\centering
\scalebox{.9}{%
\begin{tikzpicture}[x=1cm,y=1cm,font=\small,
  row/.style={rectangle,draw=figNavy,fill=white,line width=.8pt,
    minimum size=5.3pt,inner sep=0pt},
  col/.style={circle,draw=figNavy,fill=white,line width=.8pt,
    minimum size=6pt,inner sep=0pt},
  branch/.style={draw=figGray!72,line width=.65pt},
  message/.style={-{Latex[length=2mm,width=1.5mm]},line width=1pt},
  msglabel/.style={inner sep=1.5pt}
]
  \draw[figLight] (4.00,-.49) -- (4.00,5.10);

  % The two panels use identical vertices and edges.
  \foreach \panel/\initial/\boundaryColor/\parentColor/\boundaryBar/\parentBar in {
    0/0/figTeal!92!black/figOchre!82!black/\underline/\overline,
    8.00/{1/z}/figOchre!82!black/figTeal!92!black/\overline/\underline%
  } {
    \begin{scope}[xshift=\panel cm]
      \node[font=\small\bfseries,text=figNavy] at (0,5.03)
        {Boundary set to $\initial$};

      \coordinate (root) at (0,4.05);
      \coordinate (u) at (-1.52,2.08);
      \coordinate (v) at ( 1.52,2.08);
      \coordinate (u1) at (-2.93,.16);
      \coordinate (u2) at (-1.01,.16);
      \coordinate (v1) at ( 1.01,.16);
      \coordinate (v2) at ( 2.93,.16);

      % Each arrow is parallel to its edge, offset upward by 3.4 points.
      \foreach \child/\parent/\rotation in {
        u1/u/90, u2/u/-90, v1/v/90, v2/v/-90%
      } {
        \draw[branch] (\child) -- (\parent);
        \draw[message,draw=\boundaryColor]
          ($(\child)!.27!(\parent)!3.4pt!\rotation:(\parent)$) --
          ($(\child)!.72!(\parent)!3.4pt!\rotation:(\parent)$);
      }

      \foreach \child/\rotation/\labelSide in {
        u/90/above left, v/-90/above right%
      } {
        \draw[branch] (\child) -- (root);
        \draw[message,draw=\parentColor]
          ($(\child)!.27!(root)!3.4pt!\rotation:(root)$) --
          ($(\child)!.72!(root)!3.4pt!\rotation:(root)$);
        \node[msglabel,text=\parentColor,\labelSide=3pt]
          at ($(\child)!.50!(root)!3.4pt!\rotation:(root)$)
          {$\parentBar{\cH}_{\child}^{(1)}$};
        \node[col] at (\child) {};
      }

      \foreach \leaf/\vertex in {u1/u_1, u2/u_2, v1/v_1, v2/v_2} {
        \node[row] at (\leaf) {};
        \node[msglabel,below=5pt,text=\boundaryColor] at (\leaf)
          {$\boundaryBar{\cH}_{\vertex}^{(0)}=\initial$};
      }
      \node[row] at (root) {};
      \node[msglabel,above=5pt,text=\boundaryColor] at (root)
        {$\boundaryBar{\cH}_{\varnothing}^{(2)}$};
    \end{scope}
  }
\end{tikzpicture}%
}
\caption{Cavity messages travel upward on the same tree
with two different assignments at the leaves of depth $\el = 2$. Teal
denotes lower bounds and ochre denotes upper bounds. Each update
reverses order, so the all-zero boundary gives upper bounds at the
first update and a lower bound at the root. The all-$1/z$ boundary
gives the opposite bounds. At odd depths the two assignments
exchange their roles at the root.}
\label{fig:tree-messages}
\end{figure}
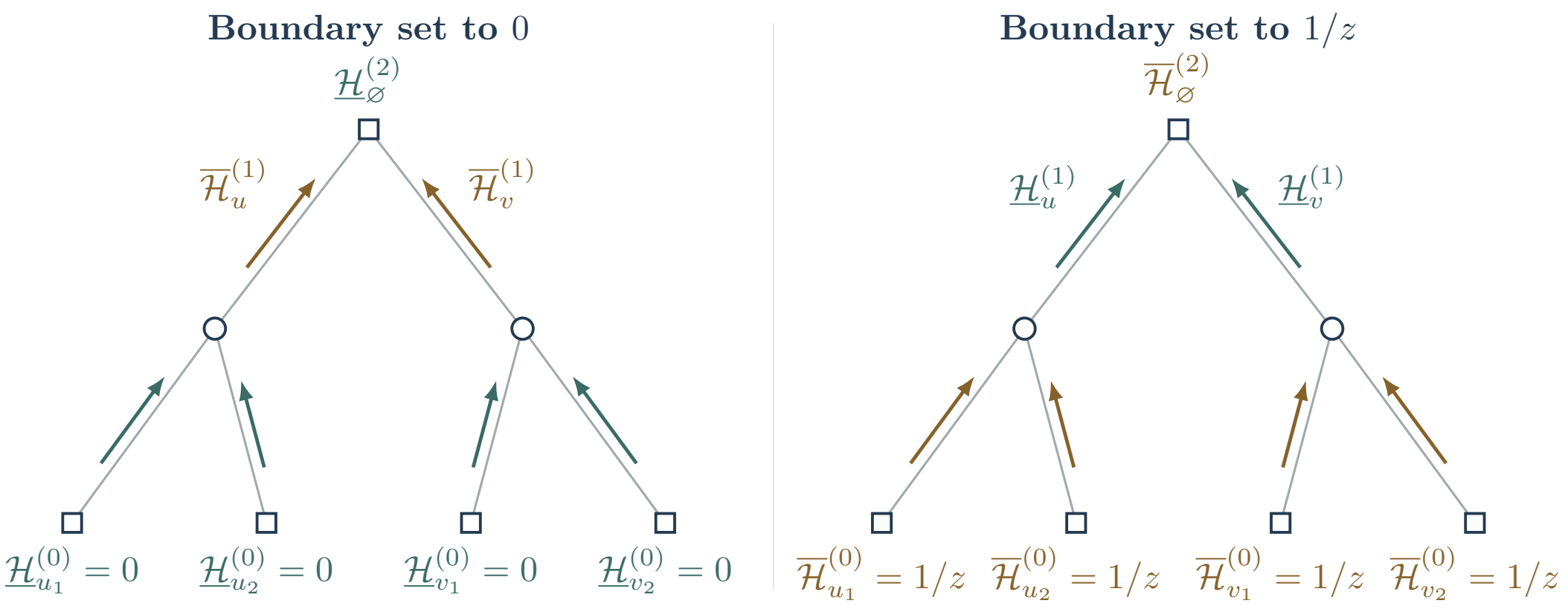

We begin with the following lemma which proves a recursion similar to the one claimed in Theorem \ref{thm:main}, only with an included monomer parameter $z > 0$. An illustration of message passing technique is shown in Figure \ref{fig:tree-messages}.

\begin{lem}[Uniqueness with monomers]\label{lem:positive}
Let $\xi = (\xi_i)_{i\geq 1}$ be a Poisson point process with intensity $\nu_d$. For every~$z>0$, the distributional equation
\begin{equation}\label{eq:positive}
 \cH_z\stackrel{\mathrm d}=\brac{z+\sum_{i\geq 1} w(\xi_i)\cH_{z,i}}^{-1}
\end{equation}
has a unique solution in~$[0,1/z]$, where $(\cH_{z,i})_{i\geq 1}$ is an independent sequence of random variables distributed as $\cH_z$ and independent of $\xi$. In particular, for any length cutoff $R > 0$, the same conclusions hold after restricting the sum to all $\xi_i$ such that $\xi_{i}\leq R$. 
\end{lem}

\begin{proof}
At every tree vertex~$v$, initialize lower and upper values~$\underline \cH_v^{(0)}=0$ and~$\ol \cH_v^{(0)}=1/z$. Define the following depth approximations recursively by
\begin{align*}
 \underline \cH_v^{(\el+1)}
 &=\brac{z+\sum_i w(\xi_{v,i})\ol \cH_{vi}^{(\el)}}^{-1},\\
 \ol \cH_v^{(\el+1)}
 &=\brac{z+\sum_i w(\xi_{v,i})\underline \cH_{vi}^{(\el)}}^{-1}\,,
\end{align*}
We first claim that $(\overline{\cH}_{v}^{(\el)})_{\el \geq 1}$ forms a non-increasing sequence in $\el$ (likewise, $(\underline{\cH}_v^{(\el)})_{\el \geq 1}$ a non-decreasing sequence in $\el$) and for every $v \in V(\cT)$, $\overline{\cH}_{v}^{(\el)} \geq \underline{\cH}_{v}^{(\el)}$. To see this, we use a simple proof by induction. To show the base case, for every depth-$\el$ vertex $v$, $\overline{\cH}_{v}^{(0)} = 1/z \geq 0 = \underline{\cH}_{v}^{(0)}$. Furthermore, by definition it is clear that $\underline{\cH}_{vi}^{(\el)} \geq \underline{\cH}_{vi}^{(\el-1)}$ implies that $\overline{\cH}_v^{(\el)} \geq \overline{\cH}_v^{(\el + 1)}$ (a similar argument implies that $\underline{\cH}_v^{(\el)} \leq \underline{\cH}_{v}^{(\el + 1)}$) and clearly $\overline{\cH}_{v}^{(\el-1)} \geq \underline{\cH}_{v}^{(\el-1)}$ implies $0 < \underline{\cH}_{v}^{(\el)} \leq \overline{\cH}_{v}^{(\el)} \leq 1/z$ for $\el \geq 1$. For the root vertex $v = \varnothing$, define the upper and lower limits of the messages as
$$
    U = \liminf_{\el \to \infty} \overline{\cH}_{\varnothing}^{(\el)}, \quad L = \limsup_{\el \to \infty} \underline{\cH}_{\varnothing}^{(\el)}.
$$
Clearly $0\leq L \leq U\leq 1/z$. Similarly define $L_i$ and $U_i$ to be the lower and upper limits of the $i$-th child of the root $\varnothing$, respectively. Notice that $(L_i,U_i)_{i\geq 1}$ are $\iid$ copies of $(L,U)$ independent of the root weights, since they are computed on disjoint subtrees in the PWIT. By summability of the outgoing weights and dominated convergence, $L$ and $U$ satisfy the joint recursions
\begin{equation}\label{eq: new rec}
    L = (z + \sum_{i} w(\xi_{\varnothing,i})U_i)^{-1}, \quad U = (z + \sum_{i} w(\xi_{\varnothing,i})L_i)^{-1}. 
\end{equation}
To show that $L = U$, we use Mecke's equation in deletion form as in \eqref{eq:mecke} with independent marks $Y_i = (U_i, L_i)$. To see this, define $L_{-i}$ be the recursive definition of the lower value in \eqref{eq: new rec} with the $i$-th vertex message omitted (so $L_{-i}$ is a function of the retained marked configuration; let $U_{-i}$ be defined similarly). Then
$$
    \frac{1}{L_{-i}} = \frac{1}{L} - w(\xi_{\varnothing,i}) U_i \implies L = \frac{L_{-i}}{1 + w(\xi_{\varnothing,i})U_iL_{-i}}.
$$
and therefore, applying Mecke's equation to the marked point process $\sum_{i} \delta_{(\xi_{\varnothing,i}, U_i, L_i)}$,
$$
1 - z\E{L} = \sum_{i} \E{w(\xi_{\varnothing,i})U_i L} = \sum_{i} \E{\frac{w(\xi_{\varnothing,i}) U_i L_{-i}}{1 + w(\xi_{\varnothing,i}) U_i L_{-i}}} = \int_0^{\infty}
\E{\frac{w(x) U^\prime L}{1 + w(x)U^\prime L}}\,\nu_d(\mathrm{d}x),$$
where $(L^\prime, U^\prime)$ is a fresh independent copy of $(L,U)$. Since the argument for $1 - z\E{U}$ is symmetric due to exchangeability, we have that $\E{U} = \E{L}$. Therefore, $U = L$ almost surely and hence $\cH_z$ is well defined. 

To show uniqueness, we show that any solution law is squeezed by these bounds. Suppose $\rho$ is a solution to $\cH_z$ on $[0,1/z]$. Starting the recursion at depth $\el$ with iid boundary messages from $\rho$, run the recursion upward as before. Since $\rho$ satisfies the same recursion as $\cH_z$, the law $\rho$ occurs at every level, including the root. However, the root is squeezed between the two depth-$\el$ bounds $\underline{\cH}_\varnothing^{(\el)}$ and $\overline{\cH}_{\varnothing}^{(\el)}$, which implies that $\cL(\cH_z) = \rho$.

To show the equivalent for the truncated tree $\cT_{R}$, apply the same message recursion only with the restriction that the sum is taken over every $i$ such that $\xi_{v,i} \leq R$. A similar argument shows that there exists a unique solution to the distributional equation $$\cH_{z,R} \overset{\mathrm{d}}{=}(z + \sum_{i: \xi_{i}\leq R}w(\xi_i)\cH_{z,R,i})^{-1}, \quad \cH_{z,R,i} \overset{\iid}{\sim} \cH_{z,R}.$$
Moreover $\lim_{R \to \infty} \cH_{z,R} = \cH_{z}$ almost surely. Indeed, run both processes on the same underlying Poisson tree $\cT$. It is clear by induction that at every fixed depth $\el$, the upper and lower messages at truncation $R$ converge to their full counterparts as $R\to\infty$. Taking $\cH_{z,R}$ to be the lower boundary limit, it is squeezed between the two cutoff bounds at every depth and hence
$$
    \underline{\cH}_\varnothing^{(\el)} \leq \liminf_{R \to \infty} \cH_{z,R} \leq \limsup_{R \to \infty}\cH_{z,R} \leq \overline{\cH}_{\varnothing}^{(\el)}.
$$
Taking $\el \to \infty$ gives the almost-sure convergence.
\end{proof}

\ssec{From neighborhoods to pressure}\label{dim:softsec}
We now give our first main result, which recovers the limiting pressure $p(\la, z)$ in terms of all solutions $(\cH_z)_{z>0}$ to the recursive distributional equation in Lemma \ref{lem:positive}. 
% Write $F(s)=\Prob{C\le s}$ and $\eta_n=n\cL(x_{11})$. The latter is the expected counting measure of lengths incident to one vertex.
% Assumption~\eqref{dim:F} gives~$\eta_n([0,x])\to x^d$. Near zero it also bounds~$F(s)$ by a constant times~$s^d$; away from zero use~$F(s)\le1$ and increase the constant. Thus there are finite constants~$C_0,C_1$, depending only on the cost law and~$d$, such that
% \begin{equation*}
%  F(s)\le C_0s^d\quad(s\ge0),
%  \qquad \eta_n([0,x])\le C_1x^d\quad(x\ge0)\,.
% \end{equation*}
% Stieltjes integration by parts, which applies without a density, gives for every~$b>0$,
% \begin{equation*}
%  \int_{(R,\infty)}e^{-bx}\,\eta_n(\dd x)
%  \le C_1b\int_R^\infty x^d e^{-bx}\dd x\,.
% \end{equation*}
% Thus long edges make a uniformly small contribution whenever their weights are bounded by an exponential in their length. This estimate also applies with an additional factor~$x$.

\begin{prop}[Pressure with monomers]\label{dim:soft}\label{prop:positive}
For fixed~$z>0$, the pressures~$p_n(\la, z)$ converge almost surely and in~$L^1$ to
\begin{equation}\label{dim:activity}\eqlabelalias{eq:soft}
 p(\la, z)=2\log z+
 2\int_z^\infty\bigl(1-u\E{\cH_u}\bigr)\frac{\dd u}{u}\,.
\end{equation}
\end{prop}

\begin{proof}
\emph{Local monomer probabilities.} Fix~$R>0$. We first show that $G_{n,R}$, rooted at a uniformly chosen
vertex, converges locally weakly to the truncated PWIT~$\cT_R$.
In other words, for each fixed radius~$\el$, the rooted ball and its
edge lengths converge in distribution to those of~$\cT_R$. This is the
bipartite version of W\"astlund's local convergence argument
\cite[Proposition~2.2 and Section~3.6]{Wastlund}.

To sketch the argument, explore the ball breadth first, revealing edges
only to undiscovered vertices in the opposite part. By~\eqref{dim:F},
the retention probability $q_{n,R}=\Prob{x_{1,1}\le R}$ satisfies
$nq_{n,R}\to R^d$. Each explored vertex therefore has asymptotically
$\operatorname{Poisson}(R^d)$ new neighbors, and the retained lengths
have limiting distribution function $(t/R)^d$ for $0\le t\le R$.
The number of vertices discovered at fixed depth is tight by counting
paths from the root. With at most~$J$ vertices discovered, the probability
of an extra retained edge among them is at most~$J^2q_{n,R}$.
Letting $n\to\infty$ and then $J\to\infty$ gives the claimed tree limit.
Starting from two prescribed vertices in opposite parts with their
connecting edge deleted gives two independent copies in the limit by
the same argument. Throughout, we write~$p_{n,R}(\la,z)=n^{-1}\log\cZ_{G_{n,R}}(z)$ for the cutoff pressure.

Next, we claim that for $v$ uniformly picked from the vertices of $G_{n,R}$,
\begin{equation}
    \lim_{n\to \infty} \E{\mu_{G_{n,R},z}(\text{$v$ unmatched})} = z\E{\cH_{z,R}}.
\end{equation}
For any finite graph, note that $\mu_{G,z}(\text{$v$ is unmatched}) = z\cH_{G,z}(v)$. For $(G_{n,R},v)$, let $L_{n,\el,R}$ and $U_{n,\el,R}$ denote the lower and upper bounds at the root $v$ computed from the radius-$\el$ ball in $G_{n,R}$ if the ball is a tree exactly as in Lemma \ref{lem:positive} (otherwise, we set $L_{n,\el,R} = 0$ and $U_{n,\el,R} = 1/z$). Since for every true child message $\cH_{i}$ such that $L_i \leq \cH_i \leq U_i$, we have that 
$$
    \frac{1}{z + \sum_{i}w_i U_i} \leq \frac{1}{z + \sum_{i}w_i \cH_i} \leq \frac{1}{z + \sum_{i}w_i L_i},
$$
the inequality propagates all the way up to the root $v$. Hence, $L_{n,\el, R} \leq \cH_{G_{n,R},z}(v) \leq U_{n,\el,R}$. Similarly, define $L_{\el,R}$ and $U_{\el,R}$ to be the bounds obtained by the same computation for the radius-$\el$ ball inside $\cT_R$. Note that for fixed $\el, R, z$, the root messages $(L,U)$ can be written as bounded, continuous functions of the radius-$\el$ neighborhoods, so by the definition of local weak convergence \cite{sznitman_objective_2004},
$$
\E{L_{n,\el,R}} \longrightarrow \E{L_{\el,R}}, \quad \E{U_{n,\el, R}} \longrightarrow \E{U_{\el, R}}.
$$
Therefore we have that 
$$
    \E{L_{\el,R}}\leq \liminf_{n\to\infty} \E{\cH_{G_{n,R},z}(v)} \leq \limsup_{n\to\infty} \E{\cH_{G_{n,R},z}(v)} \leq \E{U_{\el,R}}.
$$
By Lemma \ref{lem:positive} and the bounded convergence theorem, both sides of the inequality converge to $\E{\cH_{z,R}}$ and multiplying by $z$ proves the claim. Moreover, we also have that $\cH_{z,R}$ converges to $\cH_{z}$ almost surely as $R \to \infty$.

\emph{Integrating in the monomer factor.} For a variable monomer factor~$u>0$, define for a spanning subgraph $G$ of $G_n$
\begin{equation*}
 \ol \cH_G(u)=\frac1{2n}\sum_{v\in V(G)}\cH_{G,u}(v)\,.
\end{equation*}
Write $p_G(\la, u) = \frac{1}{n}\log \cZ_G(u)$ for the pressure. Since the probability that $v$ is unmatched is $u\cH_{G,u}(v)$,
$$
    1- u\overline{\cH}_G(u) = \frac{1}{n}\sE{|M|}{\mu_{G,u}}.
$$
Applying the differential operator $u \dv{}{u}$ to $p_G(\la, u)$,
$$
    u \dv{}{u}p_G(\la, u) = \frac{1}{n}\sE{2n - 2|M|}{\mu_{G,u}} = 2u\overline{\cH}_{G}(u).
$$
Moreover, as $u \to \infty,$ $\frac{1}{n}\log(u^{-2n}\cZ_G(u)) \to 0$ due to the empty matching contribution. Integrating $\dv{}{u}{(p_G(\la, u) - 2\log u )}$ from $z$ to $T > z$ and sending $T \to \infty$ gives
\begin{equation}
    p_G(\la, z) - 2\log z = 2\int_z^\infty (1 - u\overline{\cH}_G(u))\frac{\mathrm{d}u}{u}.
\end{equation}
To control the integral uniformly in $u$, note that the probability of an edge $e = \sbrac{s,t}$ appearing in a matching sampled from $\mu_{G,u}$ is simply ${w_e \cZ_{G-\sbrac{s,t}}(u)}/{\cZ_{G}(u)}$, so we may rewrite $1 - u\overline{\cH}_G(u)$ much more compactly as
\[
\begin{aligned}
 0 \leq 1 - u\overline{\cH}_{G}(u)
 &= \frac{1}{n}\sum_{e \in E(G)}\sE{\indic_{\sbrac{e \in M}}}{\mu_{G,u}} \\
 &= \frac{1}{n}\sum_{e=\sbrac{s,t}\in E(G)}\frac{w_e\cZ_{G-\sbrac{s,t}}(u)}{\cZ_{G}(u)}
 \leq \frac{1}{n u^2}\sum_{e \in E(G)} w_e,
\end{aligned}
\]
where we used the fact that $\cZ_G(u) \geq u^2 \cZ_{G-\sbrac{s,t}}(u)$. Focusing on the subgraph $G_{n,R}$ of $G_n$, we may rewrite in terms of the measure $\eta_n = n\mathcal{L}(x_{1,1})$ and take expectations, giving
\begin{equation}\label{eq: almost}
 0\leq \E{1 - u\overline{\cH}_{G_{n,R}}(u) } \leq \frac{1}{nu^2}\E{\sum_{i,j=1}^n e^{-\la x_{i,j}}\indic_{\sbrac{x_{i,j}\leq R}}}=u^{-2}\int_{[0,R]} e^{-\la x}\,\eta_n(\mathrm{d}x).
\end{equation}
Note that from the assumption $F(s) \sim \kappa s^d$ for $s \downarrow 0$, there exists some $C_0$ such that $F(s) \leq C_0 s^d$ near $s = 0$. Since $F(s)\le1$, increasing $C_0$ makes this bound valid for all $s\ge0$, so for every $n\ge1$ and $t\ge0$,
$$
    \eta_n([0,t]) = n F({t}/{(\kappa n)^{1/d}}) \leq \frac{C_0t^d}{\kappa}.
$$
Since $e^{-\la x} = \int_x^\infty \la e^{-\la t}\,dt$, by Tonelli's theorem, we have that
$$
    \int_0^\infty e^{-\la x}\eta_n(\mathrm{d}x) = \la \int_0^\infty e^{-\la t}\eta_n([0,t]) \mathrm{d}t \leq \frac{C_{\kappa,d}}{\la^d} < \infty
$$
The tail above $R$ is bounded by $e^{-\la R/2}\int_0^\infty e^{-\la x/2}\eta_n(\dd x)$, which tends to zero uniformly in $n$ by the same bound with $\la/2$. Thus, the inequality in \eqref{eq: almost} is bounded from above by $\frac{C}{u^2}$ for some constant $C$ independent of $n,R,u$. Therefore, by an application of dominated convergence, we have the following limit for fixed $R$ as $n\to\infty$
\begin{equation}
    p_R(\la, z):=\lim_{n \to \infty} \E{p_{n,R}(\la, z)} = 2 \log z + 2 \int_z^\infty (1 - u\E{\cH_{u,R}})\frac{\mathrm{d}u}{u}.
\end{equation}
A second application of dominated convergence theorem with respect to $R$ yields
\begin{equation}\label{eq:final}
    p(\la, z) := \lim_{R \to \infty} p_R(\la, z) = 2\log z + 2\int_z^\infty (1 - u\E{\cH_u})\frac{\mathrm{d}u}{u}.
\end{equation}
It remains to compare the full pressure $p_n(\la, z)$ to the cutoff pressures $p_{n,R}(\la, z)$. If $e = \sbrac{s,t}$ is a particular edge absent in $G$, then the gain in the partition function $\cZ_{G + e}(z)$ is exactly $w_e \cZ_{G-\sbrac{s,t}}(z)$. This implies that
$$
    \log \cZ_{G + e}(z) - \log \cZ_G(z) = \log(1 + \frac{w_e \cZ_{G-\sbrac{s,t}}(z)}{\cZ_G(z)}) \leq \frac{w_e}{z^2}.
$$
If $\sbrac{e_1,\dots,e_k}$ are the absent edges from $G_n$ in $G_{n,R}$, applying the above argument iteratively, dividing by $n$, and taking expectations gives
$$
    0\leq \E{|p_n(\la, z) - p_{n,R}(\la, z)|} \leq z^{-2}\int_{(R,\infty)} e^{-\la x}\,d\eta_n(x).
$$
Consequently, for every fixed $R$, we get that
$$
\limsup_{n\to\infty} |\E{p_n(\la, z)} - p(\la, z)| \leq z^{-2}\sup_{n\ge1}\int_{(R,\infty)} e^{-\la x}\,d\eta_n(x) + |p_R(\la, z) - p(\la, z)|.
$$
By the uniformity of the tail and \eqref{eq:final}, we prove the convergence of the expected pressures.

\emph{Concentration by independent rows.}
Fix $z>0$. Replacing one cost row changes only the weights of
edges incident to its row vertex $v$. Let $G$ and $G'$ be the
weighted graphs before and after the replacement. Deleting $v$
gives the same weighted graph in both cases, so write
$D=\cZ_{G-v}(z)=\mathcal Z_{G'-v}(z)$.

The vertex-deletion identity \eqref{eq:delete} splits the partition
function according to whether $v$ is unmatched or matched to a
neighbor $a$:
\[
 \mathcal Z_G(z)
 = zD+\sum_{a\sim v}w_{va}\mathcal Z_{G-\sbrac{v,a}}(z).
\]
For each $a$, the matchings of $G-v$ that leave $a$ unmatched
contribute $z\mathcal Z_{G-\sbrac{v,a}}(z)$ to $D$. Hence
$\mathcal Z_{G-\sbrac{v,a}}(z)\le D/z$. Since $v$ has at most $n$ neighbors
and $w_{va}\le1$, it follows that
\[
 zD\le \mathcal Z_G(z)\le (z+n/z)D.
\]
The same bounds hold for $\mathcal Z_{G'}(z)$ with the same $D$.
Taking logarithms and dividing by $n$, we see that replacing a
single row changes the pressure by at most
\[
 c_n(z)=\frac1n\log\!\left(\frac{z+n/z}{z}\right)
       =\frac1n\log(1+n/z^2).
\]
The cost rows are independent, so McDiarmid's inequality
\cite{McDiarmid} gives, for every $t>0$,
\begin{equation}\label{eq:soft-concentration}
 \mathbb P\bigl(|p_n(\la, z)-\E{p_n(\la, z)}|>t\bigr)
 \le 2\exp\!\left(-\frac{2t^2}{n c_n(z)^2}\right)
 = 2\exp\!\left(-\frac{2nt^2}{\log^2(1+n/z^2)}\right).
\end{equation}
For fixed $z,t>0$, the bound is at most $2e^{-2t^2\sqrt n}$ for
all sufficiently large $n$, and is therefore summable. Applying
the first Borel--Cantelli lemma at each threshold $t=1/k$,
$k\ge1$, shows that $|p_n(\la, z)-\E{p_n(\la, z)}|\to0$ almost surely. Integrating \eqref{eq:soft-concentration} over $t>0$ also gives
\[
 \E{|p_n(\la, z)-\E{p_n(\la, z)}|}
 \le 2\int_0^\infty
      \exp\!\left(-\frac{2t^2}{n c_n(z)^2}\right)\,\mathrm dt
 = \sqrt{\frac{\pi}{2}}\,
   \frac{\log(1+n/z^2)}{\sqrt n}
 \longrightarrow 0.
\]
Since $\E{p_n(\la, z)}\to p(\la, z)$, the triangle inequality gives
$p_n(\la, z)\to p(\la, z)$ almost surely and in $L^1$.
\end{proof}

\ssec{Final limiting formula}\label{subsec:bethe}
In this section, we rewrite the formula of the limiting pressure at a fixed monomer activity $z > 0$ derived in Proposition \ref{prop:positive} using only the cavity message \textit{at that $z$} (without integrating over all activities $u\geq z$ required in the previous formula for $p(\la, z)$). This is similar in spirit to the final form of the limiting pressure appearing in Theorem \ref{thm:main}.
% \begin{lem}\label{lem:bethe}
% For~$z>0$, let $\cH_z'$ is an independent copy of~$\cH_z$, where $\cH_z$ is the root message on the Poisson tree. Then,
% \begin{equation}\label{dim:bethe}\eqlabelalias{eq:bethe}
%  p(\la, z)=-2\E{\log \cH_z}-
% \int_0^\infty\E{\log\bigl(1+w(x)\cH_z\cH_z'\bigr)}\,\nu_d(\dd x)\,.
% \end{equation}
% \end{lem}

\begin{lem}[Vertex--edge identity]\label{lem:bethe}
Let $\cH_z$ be the root message on the Poisson tree at monomer
activity $z>0$, and let $\cH_z'$ be an independent copy. Then
\begin{equation}\label{dim:bethe}\eqlabelalias{eq:bethe}
 p(\la, z)
 =-2\E{\log\cH_z}
 -\int_0^\infty
    \E{\log\bigl(1+w(x)\cH_z\cH_z'\bigr)}\,\nu_d(\mathrm dx),
\end{equation}
where $w(x)=e^{-\lambda x}$ and
$\nu_d(\mathrm dx)=d x^{d-1}\,\mathrm dx$.
\end{lem}

\begin{proof}
Write $\Phi(z)$ for the right-hand side of \eqref{dim:bethe}.
At the root, let $w_i=w(\xi_i)$ and $S=\sum_i w_i$. Then
\[
 K:=\E{S}=\int_0^\infty w(x)\,\nu_d(\mathrm dx)<\infty.
\]
The root recursion \eqref{eq:positive} and the bound
$0<\cH_{z,i}\le 1/z$ on the child messages give
\[
 0\le -\log(z\cH_z)
 =\log\left(1+\frac{1}{z}\sum_i w_i\cH_{z,i}\right)
 \le \frac{S}{z^2}.
\]
Also,
$\log(1+w(x)\cH_z\cH_z')\le w(x)/z^2$.
Therefore, both terms of $\Phi(z)$ are finite.

Next, we construct messages for different activities $z$ and $y$ on the
same realization of the tree $\cT$. For $0<z\le y$, notice that the finite-depth
approximations in Lemma~\ref{lem:positive} at level $\el$ satisfy
\[
 \frac{z}{y}\cH_z^{(\ell)}
 \le \cH_y^{(\ell)}
 \le \frac{y}{z}\cH_z^{(\ell)},
\]
where we have assumed that at the boundary we take the values $0$ and $1/z$ at activity $z$ (respectively, $0$ and $1/y$ at activity $y$) and we make the same choice of either lower or upper approximation
for both activities.
By Lemma~\ref{lem:positive}, the depth limits as $\el \to \infty$ exist almost surely at every vertex simultaneously for a countable dense set of rational activities $\mathbb{Q}_{> 0}$. On this dense set, for general $0 < z\leq y$, taking an approximating sequence in $\mathbb{Q}_{> 0}$ to define $$\log \cH_v(z) := \lim_{\substack{q_n \to z\\q_n \in \mathbb{Q}_{> 0}}}\log \cH_{v}(q_n)$$ and logarithms gives almost surely
\begin{equation}\label{eq: comparis}
    \bigl|\log\cH_y-\log\cH_z\bigr|
 \le \log(y/z)\le \frac{y-z}{z}.
\end{equation}
In particular, the extended log messages are $1/\alpha$-Lipschitz in $z$ on every interval $[\alpha,\beta]\subset(0,\infty)$. Using an approximation argument, both depth approximations converge to $\cH_v(z)$ at every positive activity on the same probability-one event. By mean value theorem, applying the above inequality, we may conclude that $\E{\cH_z}$ is locally Lipschitz, hence continuous everywhere and differentiable almost everywhere. Differentiating the activity formula \eqref{dim:activity} then gives
\[
 \partial_z p(\la,z) = \frac{2}{z}-\frac{2}{z} (1-z\E{\cH_z}) = 2\E{\cH_{z}}.
\]
Write
\[
 D_z=\frac{\mathrm d}{\mathrm dz}\log\cH_z,
 \qquad
 D_{z,i}=\frac{\mathrm d}{\mathrm dz}\log\cH_{z,i}.
\]
Note that \eqref{eq: comparis} implies $|D_z|\le 1/z$ and
$|D_{z,i}|\le 1/z$ wherever the derivatives exist. By local
Lipschitz continuity and Fubini's theorem, for almost every $z$ the derivatives exist almost surely at every vertex.
On $[\alpha,\beta]$, the messages themselves are Lipschitz
with constant at most $1/\alpha^2$, since their logarithms are
$1/\alpha$-Lipschitz and their values are bounded by $1/\alpha$.
Thus the difference quotients in the weighted child sum are
dominated by the summable sequence $(w_i/\alpha^2)_i$.
Differentiating the logarithm of the root recursion gives
\[
 -D_z
 =\frac{1+\sum_i w_i\cH_{z,i}D_{z,i}}
        {z+\sum_i w_i\cH_{z,i}}
 =\cH_z+\sum_i w_i\cH_{z,i}\cH_z D_{z,i}.
\]
To take the expectation of the sum, let
\[
 A_{-i}=\left(z+\sum_{j\ne i} w_j\cH_{z,j}\right)^{-1}
\]
be the root message with edge $i$ omitted. Reinserting that edge
adds $w_i\cH_{z,i}$ to the inverse message. Hence by a similar sequence of steps as in Lemma \ref{lem:positive},
\[
 \cH_z=\frac{A_{-i}}{1+w_iA_{-i}\cH_{z,i}},
 \qquad
 w_i\cH_{z,i}\cH_z
 =\frac{w_iA_{-i}\cH_{z,i}}{1+w_iA_{-i}\cH_{z,i}}.
\]
We now apply the deletion form of Mecke's identity \eqref{eq:mecke}
with marks given by $(\cH_{z,i},D_{z,i})$. In the insertion integral,
the remaining point process supplies the first root message
$\cH_z$ and the inserted point carries an independent mark
$(\cH_z',D_z')$ from a second tree. Writing
\[
 a_x(z)=\frac{w(x)\cH_z\cH_z'}{1+w(x)\cH_z\cH_z'},
\]
we obtain
\[
 -\E{D_z}
 =\E{\cH_z}
  +\int_0^\infty \E{a_x(z)D_z'}\,\nu_d(\mathrm dx),
\]
where we have justified applying Mecke's formula since $a_x(z) \leq w(x)/z^2$ and $|D_z'|\leq 1/z$ implies $\int_0^\infty \E{|a_x(z) D_z'|}\,d\nu_d(x)\leq z^{-3}\int_0^\infty w(x)\,d\nu_d(x) < \infty$.

For the edge integral in $\Phi(z)$, use two independent copies
of the common-tree construction. The chain rule gives
\[
 \frac{\mathrm d}{\mathrm dz}
 \log\bigl(1+w(x)\cH_z\cH_z'\bigr)
 =a_x(z)(D_z+D_z').
\]
On $[\alpha,\beta]$, its absolute value is bounded by
$2w(x)/\alpha^3$. This bound and the local Lipschitz bound on
$\log\cH_z$ justify differentiating under the integrals and show
that $\Phi$ is locally absolutely continuous. Since exchanging
the two trees preserves $a_x(z)$, the expectations of
$a_x(z)D_z$ and $a_x(z)D_z'$ are equal. Consequently, for
almost every $z$,
\[
 \Phi'(z)
 =-2\E{D_z}
   -2\int_0^\infty\E{a_x(z)D_z'}\,\nu_d(\mathrm dx)
 =2\E{\cH_z}
 =\partial_z p(\la,z).
\]
Thus $\Phi(z)-p(\la,z)$ is constant. To show that the constant difference is zero, notice that the initial bound on $\Phi(z)$ shows that 
$|\Phi(z)-2\log z|\lesssim \frac{K}{z^2}\to0$ as $z\to\infty$. Similarly, the root recursion at activity $u$ also gives
\[
 0\le 1-u\cH_u
 =\cH_u\sum_i w_i\cH_{u,i}\le \frac{S}{u^2}.
\]
Substituting this into \eqref{dim:activity} yields
\[
 0\le p(\la, z)-2\log z
 \le 2K\int_z^\infty \frac{\mathrm du}{u^3}
 =\frac{K}{z^2}.
\]
Therefore $\Phi(z)-p(\la, z)\to0$ as $z\to\infty$, and therefore $\Phi(z) = p(\la, z)$ as desired.
\end{proof}
\section{Weighted completion}\label{sec:completion}
In this section, our goal will be to rigorously prove $\lim_{z \downarrow 0} p(\la, z)$ is equal to the almost-sure limit  $\lim_{n\to\infty}p_n(\la)$. To accomplish this, our technique uses a  \textit{weighted completion} argument to complete partial matchings along short, low-cost augmenting paths, then bound the total weight of all partial matchings which complete a given perfect matching. To recover perfect assignments, we choose a particular surjection (via augmenting paths) from partial to perfect matchings $\mathsf{P}(M)$ and regroup the weights by that choice: for each $k \in [n]$
\begin{equation}\label{eq:regroup}
 \cZ_{n,k}=\sum_{P\text{ perfect}}\weight_\la(P)
 \sum_{\substack{|M|=n-k\\\mathsf{P}(M)=P}}
 \frac{\weight_\la(M)}{\weight_\la(P)}\,.
\end{equation}
The task is to bound the inner sum of \eqref{eq:regroup} uniformly for every perfect matching $P \in \mathscr{M}_n$. We first control the edges changed by completion, then sum the weight ratios over all possible changed rows. Importantly, the random events used below depend only on the cost array $\mathbf C$ and hold for every partial matching simultaneously. The alternating-path construction follows the type of completion arguments used in ~\cite[Section~3.3, Lemmas~3.5--3.7]{Wastlund}. 

\ssec{Completion along short alternating paths}

\begin{figure}[!htbp]
\centering
\begin{tikzpicture}[x=1cm,y=1cm,font=\small,>=Latex,
 row/.style={rectangle,draw=figNavy,fill=figNavy,minimum size=5.3pt,inner sep=0pt},
 col/.style={circle,draw=figTeal,fill=figTeal,minimum size=5.8pt,inner sep=0pt},
 old/.style={draw=figNavy,line width=1.0pt},
 cheap/.style={draw=figOchre,line width=1.0pt,dashed}]
 \node[anchor=south] at (2.925,2.25) {Before augmentation};
 \node[anchor=south] at (11.025,2.25) {After augmentation};
 \node[anchor=east] at (.35,1.4) {$R_n$};
 \node[anchor=east] at (.35,0) {$L_n$};
 \foreach \k in {0,1,2}{
   \draw[old] ({.9+1.35*(\k+1)},1.4)--(.9+1.35*\k,0);
 }
 \foreach \k in {0,1,2,3}{
   \draw[cheap] (.9+1.35*\k,1.4)--(.9+1.35*\k,0);
   \draw[draw=figOchre,line width=1.05pt] (9+1.35*\k,1.4)--(9+1.35*\k,0);
   \node[row] at (.9+1.35*\k,1.4) {};
   \node[col] at (.9+1.35*\k,0) {};
   \node[row] at (9+1.35*\k,1.4) {};
   \node[col] at (9+1.35*\k,0) {};
 }
 \node[row,fill=white,line width=.9pt,minimum size=7pt] at (.9,1.4) {};
 \node[col,fill=white,line width=.9pt,minimum size=8pt] at (4.95,0) {};
 \node[above=5pt,font=\footnotesize] at (.9,1.4) {free};
 \node[below=5pt,font=\footnotesize] at (4.95,0) {free};
 \draw[-{Latex[length=2mm]},figGray,line width=.65pt] (5.85,.7)--(7.45,.7);
 \node[above=4pt,font=\footnotesize] at (6.65,.7) {flip};
 \node[below=4pt,font=\footnotesize,align=center] at (6.65,.7) {$+1$ edge};
\end{tikzpicture}
\caption{An illustrative alternating path: squares lie in~$R_n$, circles in~$L_n$. Before the flip, matching edges are solid and cheap unused edges are dashed. The flip removes three edges and inserts four, covering the two hollow free endpoints while keeping every interior vertex matched. Only the path is shown.}\label{fig:augment}
\end{figure}
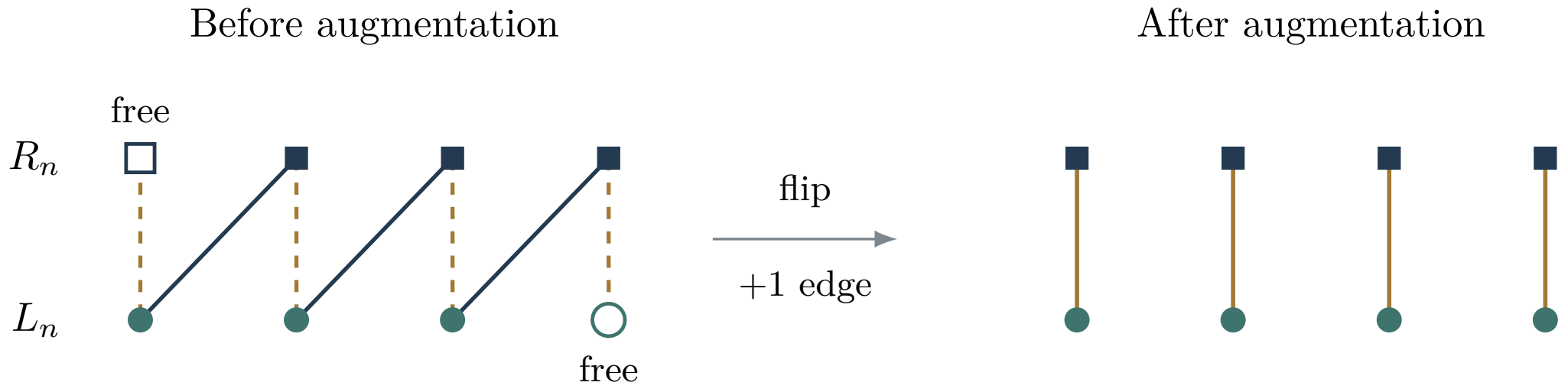

Before we begin, we need to go over the necessary formalism on weighted completions. For a matching~$M$, we define its \emph{deficit} as ~$n-|M|$, the number of unmatched vertices in each class of vertices. A vertex is \emph{free} if no edge in a matching is incident to it. We define an \emph{augmenting path} in $G_n$ relative to $M$ to be a path without repeated vertices from a free row to a free column vertex, alternating edges outside and inside the matching, and beginning and ending outside it. Flipping it adds one matched edge (Figure~\ref{fig:augment}). At deficit~$j$, a length cutoff of order~$\log^{1/d}(en/j)$ permits~$O(\log(en/j))$ doubling rounds to reach a fixed fraction of each class. One augmentation therefore inserts~$O(\log(en/j))$ edges of total length~$O(\log^{1+1/d}(en/j))$. Summing over the deficits gives the next lemma.

Let $A > 0$ be a constant to be chosen. For integer deficits $j \in \sbrac{1,2,\dots,n}$, let $D_{j,n}$ be the bipartite graph obtained by discarding all edges of $G_n$ whose weight is greater than $r_{j,n}:= A \log^{1/d}(en/j)$ and $\pi_{j,n} = \Prob{x_e \leq r_{j,n}}$ be the retention probability at level $r_{j,n}$. For a given set of vertices $S$ in either $L_n$ or $R_n$, let $N_{D_{j,n}}(S)$ be the neighboring vertices of the opposite type joined to at least one vertex of $S$ in $D_{j,n}$. Set $h = \lfloor n/40\rfloor$.\footnote{The particular choice of $h$ is merely convenient; the same proof of Lemma \ref{lem:expander} works with any $h = \lfloor cn \rfloor$ for $0 < c < 1/10$ and suitable $A$.} We begin with a necessary expander lemma which we show holds simultaneously for all levels $j$.
\begin{lem}\label{lem:expander}
    For every $c > 0$, except on an event of probability $O(n^{-c})$, the following occurs simultaneously for every $j \in \sbrac{1,2,\dots, n}$:
    \begin{enumerate}
        \item If $S\subseteq L_n$ or $S\subseteq R_n$, then $|N_{D_{j,n}}(S)|\geq 10 |S|$ whenever $j \leq |S| \leq h$. 
        \item Any two sets of vertices from $L_n$ and $R_n$ of size $h$ must have an edge of $D_{j,n}$ going between them.
    \end{enumerate}
\end{lem}
\begin{proof}
Since~$r_{1,n}/(\ka n)^{1/d}\to0$, assumption~\eqref{dim:F} gives, uniformly in~$1\le j\le n$ for all sufficiently large~$n$,
\begin{equation*}
 \pi_{j,n}=F\bigl(r_{j,n}/(\ka n)^{1/d}\bigr)
 \ge\frac{A^d}{2n}\log(en/j)\,.
\end{equation*}
For fixed $s$ satisfying $j\le s\le h$, failure of Property $1$ in Lemma \ref{lem:expander} means that some size $|S| = s$ set of vertices in either $L_n$ or $R_n$, all of its neighbors in $N_{D_{j,n}}(S)$ lie in a set of size strictly less than $10s$. The union bound over possible choices of subsets $S$ and neighbors $N_{D_{j,n}}(S)$ gives
\begin{align}\label{eq: property1}
 2\binom ns\binom n{10s}(1-\pi_{j,n})^{s(n-10s)}
 &\le2\exp\!\left[-\brac{\frac{3A^d}{8}-11}s\log(en/s)\right].
\end{align}
Indeed, using $\binom{n}{m}\leq (en/m)^m$, $\log(\binom{n}{s}\binom{n}{10s})\leq 11s\log(en/s)$. Since $s \leq h = \lfloor n/40\rfloor$ and $j\leq s$, $n-10s \geq 3n/4$, we have that $\pi_{j,n}(n-10s)s \geq (3A^{d}/8)s\log(en/s)$. Therefore, multiplying the two factors gives \eqref{eq: property1}
as required. Summing over all choices of deficit $j$ and sizes $j\leq s \leq h$, it suffices to choose $A$ sufficiently large enough so that $3A^d/8 - 11 > c + 2$ to have the probability vanishing at the order $O(n^{-c})$. 

For the failure of Property 2 in Lemma \ref{lem:expander}, note that there are at most~$4^n$ pairs of opposite sets. For n sufficiently large ($n\ge80$), use~$h\ge n/80$ and~$\pi_{j,n}\ge A^d/(2n)$ to obtain, for each~$j$,
\begin{equation*}
 4^n(1-\pi_{j,n})^{h^2}
 \le\exp\!\left[n\log4-\frac{h^2A^d}{2n}\right]\, = \exp\!\left[n\brac{\log4-\frac{A^d}{12800}}\right]\, .
\end{equation*}
Again, choosing~$A$ sufficiently large so that~$A^d/12800>\log4+1$ and summing over all possible values of $j$ gives an exponentially small failure probability (so $O(n^{-c})$ for every $c > 0$).
\end{proof}

\begin{lem}[Cost of completion]\label{lem:augment}
Fix pseudodimension $d>0$, $c >0$, and suppose the cost $C$ has a law satisfying~\eqref{dim:F}. There is~$K_0<\infty$ depending on $d$ and $c$, such that, outside an event of probability~$O(n^{-c})$, every matching~$M_j \in \mathscr{M}(G_n)$ of size~$n-j$ for some $1\leq j \leq n$ has a deterministically chosen perfect completion~$P=\mathsf{P}(M_j)$ satisfying
\begin{equation}\label{eq:augment}
 \begin{split}
 |P\setminus M_j|&\le K_0j\log(en/j),\\
 \sum_{e\in P\setminus M_j}x_e&\le K_0j\log^{1+1/d}(en/j).
 \end{split}
\end{equation}
\end{lem}

\begin{figure}[!htbp]
  \centering
  \begin{tikzpicture}[x=1cm,y=1cm,font=\small,
    row/.style={rectangle,draw=figNavy,fill=white,line width=.8pt,
      minimum size=5.3pt,inner sep=0pt},
    col/.style={circle,draw=figTeal,fill=figTeal,line width=.8pt,
      minimum size=5.8pt,inner sep=0pt},
    root/.style={col,fill=white,minimum size=6.8pt},
    matched/.style={draw=figNavy,line width=1pt},
    cheap/.style={draw=figOchre,line width=1pt,dashed},
    round brace/.style={decorate,decoration={brace,amplitude=4pt},
      draw=figGray,line width=.55pt},
    vertex label/.style={anchor=east,fill=white,inner sep=1.5pt}]
    \foreach \y/\side in {4.32/L,3.32/R,2.32/L,1.32/R,.32/L}{
      \draw[figLight,line width=.35pt] (0,\y) -- (8.72,\y);
      \node[anchor=east] at (0,\y) {$\side_n$};
    }

    % In each round, an old outer vertex attaches a new matched pair.
    \foreach \parent/\child in {2.94/1.86,2.94/4.02,5.64/5.64,7.8/7.8}{
      \draw[cheap] (\parent,4.32) -- (\child,3.32);
      \draw[matched] (\child,3.32) -- (\child,2.32);
      \node[row] at (\child,3.32) {};
    }
    \foreach \parent/\child in {
      1.86/.78,1.86/1.86,1.86/2.94,4.02/4.02,
      5.64/5.1,5.64/6.18,7.8/7.26,7.8/8.34%
    }{
      \draw[cheap] (\parent,2.32) -- (\child,1.32);
      \draw[matched] (\child,1.32) -- (\child,.32);
      \node[row] at (\child,1.32) {};
      \node[col] at (\child,.32) {};
    }
    \foreach \x in {1.86,4.02,5.64,7.8}{\node[col] at (\x,2.32) {};}
    \foreach \x/\k in {2.94/1,5.64/2,6.54/3,7.8/j}{
      \node[root] at (\x,4.32) {};
      \node[above=3pt] at (\x,4.32) {$r_{\k}$};
    }
    \node[fill=white,inner sep=1pt] at (7.17,4.32) {$\cdots$};
    \node[vertex label] at (1.73,2.32) {$a$};
    \node[vertex label] at (.65,1.32) {$y$};
    \node[vertex label] at (.65,.32) {$z$};

    \draw[round brace] (8.94,3.88) -- (8.94,2.18);
    \draw[round brace] (8.94,1.91) -- (8.94,.18);
    \node[anchor=west,align=left] at (9.24,3.03)
      {Round 1\\$j\to 2j$ outer vertices};
    \node[anchor=west,align=left] at (9.24,1.045)
      {Round 2\\$2j\to 4j$ outer vertices};

    \draw[cheap] (1.8,-.47) -- (2.5,-.47)
      node[right=3pt] {$D_{j,n}\setminus M_j$};
    \draw[matched] (5.22,-.47) -- (5.92,-.47)
      node[right=3pt] {$M_j$};
  \end{tikzpicture}
  \caption{Two rounds of growth of the alternating forest $\cF_L$ rooted in $L_n$. Open circles are the $j$ free roots, filled circles are matched outer vertices, and squares are inner vertices in $R_n$. Each attachment follows a dashed edge in $D_{j,n}\setminus M_j$ and then a solid edge in $M_j$, adding one outer vertex. The increase in counts $j\to 2j\to 4j$ applies only when $4j\le h$.}
  \label{fig:alternating-forest-growth}
\end{figure}
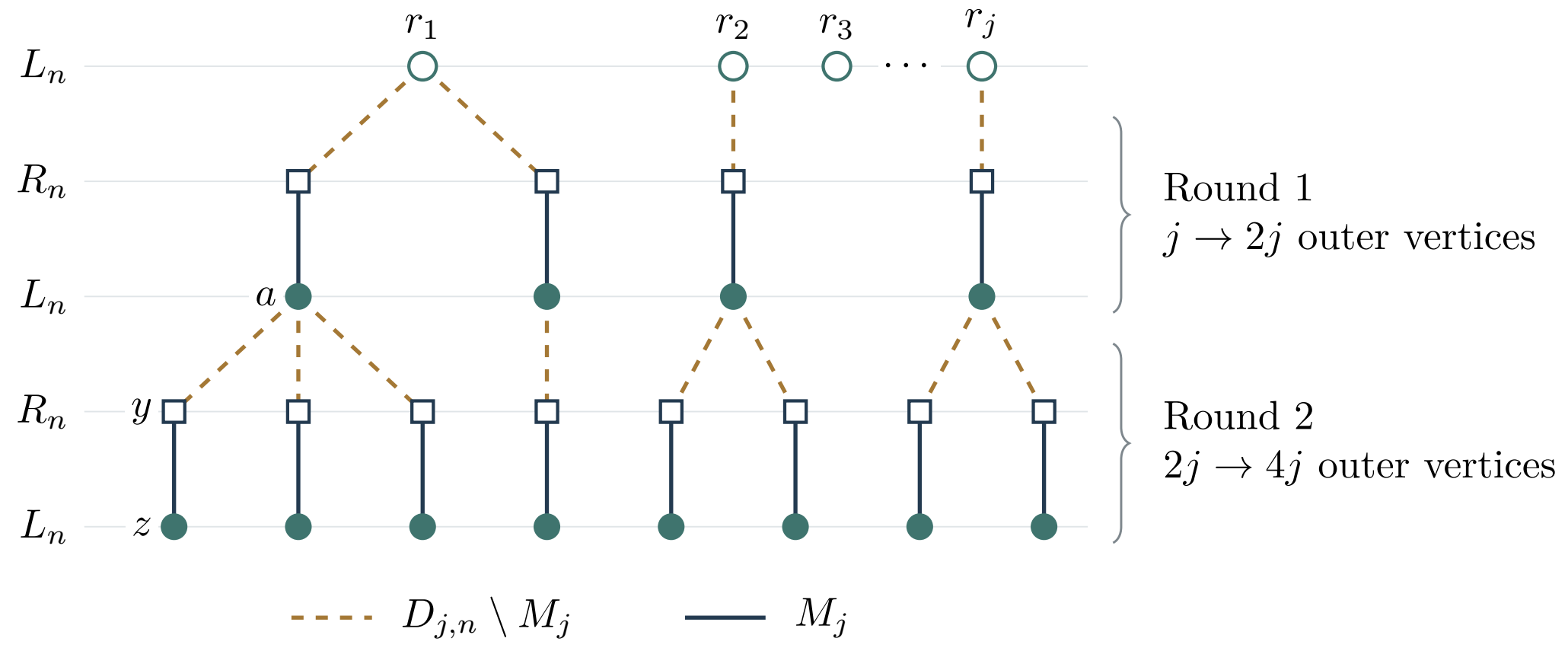

\begin{proof}
    
Suppose $1 \leq j < h$. For a fixed matching $M_j$ of size $n-j$, let $U_L \subseteq L_n$ and $U_R\subseteq R_n$ be its free vertices of size $j$. Let $\cF_L$ and $\cF_R$ be two disjoint forests started at the isolated vertices of $U_L$ and $U_R$ respectively, where the roots of each tree are the $j$ different free vertices in $U_L$ and $U_R$ respectively. Our algorithm enlarges each forest over the course of multiple rounds, attaching new vertices such that edges alternate between membership in $D_{j,n}\setminus M_j$ and $M_j$. We define the \textit{outer vertices} of the forests $\cF_L$ and $\cF_R$ to be the vertices of their tree components belonging to the same type (either in $L_n$ or $R_n$, see Figure \ref{fig:alternating-forest-growth}). We will continue over multiple rounds, adding vertices until the total number of outer vertices reaches $h$. Suppose that at the beginning of a round, both outer sets of vertices in the trees $\cF_L$ and $\cF_R$ have size $s$ where $j\leq s < h$. Set $s^\prime  = \min(2s, h)$. We first enlarge the outer set of vertices from size $s$ to size $s^\prime$ in $\cF_L$ (and later do the same procedure applied to $\cF_R$). For every edge added to the forest from $D_{j,n}\setminus M_j$, the edge from $M_j$ determines exactly $2$ new vertices added to the tree. Therefore, since there are $j$ roots and $2(s - j)$ vertices coming from edges in $M_j$, before the start of the round there are a total of $2s - j$ vertices in the forest $\cF_L$. Therefore, after one round applied to $\cF_L$, the total number of vertices increases to $2s^\prime - j + 2s - j \leq 6s$ in total. To justify the existence of such rounds, suppose $S \subseteq L_n$ is the current set of open vertices in $\cF_L$. Since $j \leq |S|\leq h$, by Property 1 in Lemma \ref{lem:expander}, we have that $|N_{D_{j,n}}(S)|\geq 10|S|$ and hence there are at least $4|S|$ new neighbors after discarding vertices in either forest. 

Take one such unused neighbor $y$ of an open vertex $a$ (see Figure \ref{fig:alternating-forest-growth}). Note that it must be matched in $M_j$ since the $j$ roots of the forests are the only unmatched vertices in $G_n$ (likewise, the corresponding vertex $z$ such that $(y,z)\in M_j$ must also be unused, because otherwise it would have been seen). Since $4|S| \geq s^\prime - s$, there are enough possible neighbors. Since we perform the expansion relative to the type of vertices in $S$, at each stage the new pairs $(y,z)\in M_j$ added are disjoint across $\cF_L$ and $\cF_R$. Applying this procedure repeatedly over $t_j := \lceil \log_2(h/j) \rceil$ rounds, we stop the procedure once $\cF_L$ and $\cF_R$ have exactly $h$ outer vertices, noting that the length from a root to any outer vertex is at most $2t_j$. By Property 2 of Lemma \ref{lem:expander}, there exists an edge in $D_{j,n}\setminus M_j$ between an outer vertex $u$ from $\cF_L$ and $v$ from $\cF_R$. If $r_u$ and $r_v$ are the free roots corresponding to $u$ and $v$ respectively, then the path $r_u \to u \to v \to r_v$ forms an augmenting path in $G_n$ of length at most $2t_j + 1$. These paths will make up the augmenting path surjection $M \mapsto \mathsf{P}(M)$. In the case that $j \geq h$, we will simply choose a cheap edge in $D_{j,n}\setminus M_j$ between any arbitrary pair of size $h$ subsets. 

To complete the proof, we now sum the contributions along augmentations. To build the perfect matching $P = \mathsf{P}(M_j)$, we fix an ordering on the vertex labels and break all ties lexicographically. Starting at deficit $j$ with $M^{(j)} = M_j$, we repeatedly apply the above construction inductively at decreasing deficits $k = j, j-1, \dots, 1$ to matchings $\{M^{(k)}\}_{1\leq k \leq j}$, where $M^{(k-1)}$ is the matching after flipping the augmenting path found by applying the above algorithm to the matching $M^{(k)}$ with deficit $k$. We define the perfect matching at the end of this process as $\mathsf{P}(M_j) = M^{(0)} = P$. Note that the expansion step using Property $1$ of Lemma \ref{lem:expander} holds across all $j \in \sbrac{1,\dots,n}$. 

At deficit $k$, let $Q_k$ be the augmenting path between two free vertices of the construction. After flipping along the augmenting path $Q_k$, we removed all edges belonging to $Q_k\cap M^{(k)}$ and inserted edges belonging to $D_{k,n}\setminus M^{(k)}$. Let $E_k = M^{(k-1)}\setminus M^{(k)}$ be the edges inserted at deficit $k$. Since $h \leq n$, the bound $2t_k + 1$ for $k < h$ together with the one-edge bound for $k \geq h$, gives a constant $C > 0$ such that
\begin{equation}\label{eq:contribution-k}
    |E_k| \leq C\log(en/k) \implies \sum_{e \in E_k} x_e \leq C A \log^{1+1/d}(en/k).
\end{equation}
Since every edge in $\mathsf{P}(M_j)\setminus M_j$ must have entered the matching at some step, we have that $\mathsf{P}(M_j)\setminus M_j \subseteq \bigcup_{k=1}^j E_k$. Note that for fixed $b > 0$, 
$$
    \sum_{j=1}^k \log^b(en/j) = O(k\log^b(ne/k)),
$$
which after applying this estimate to \eqref{eq:contribution-k} for $b = 1$ and $b = 1 + 1/d$ and summing over $1\leq k \leq j$ gives for some $K_0 > 0$,
\begin{equation}
    \abs{\mathsf{P}(M_j)\setminus M_j}\leq K_0 j\log(en/j), \quad \sum_{e \in \mathsf{P}(M_j)\setminus M_j} x_e \leq K_0 j \log^{1+1/d}(en/j),
\end{equation}
as required. 
\end{proof}

\ssec{Counting weighted preimages}

From now on, we fix the following exponent $\alpha_d$ throughout:
\begin{equation*}
 \al_d=\max\{2,1+1/d\}\,.
\end{equation*}

\begin{lem}[Weighted completion]\label{dim:completion}\label{lem:weights}
For every~$c>0$ there is~$C<\infty$, depending on~$c,d,\la$ and the cost law, such that, except on an event of probability~$O(n^{-c})$, simultaneously for~$1\le j\le n$,
\begin{equation}\label{dim:completion-bound}\eqlabelalias{eq:ratio}
 \cZ_{n,j}\le \cZ_n(\la)\exp\!\brac{Cj\log^{\al_d}(en/j)}\,.
\end{equation}
\end{lem}

\Needspace{6\baselineskip}
\begin{proof}
Suppose we are working on the event in \lemref{lem:augment} with~$K_0\ge1$ and let $\cE_0$ be such event. In addition, if we define $S_i := \sum_{e \ni i_R} w_e$ to be the total weights of edges emanating from $i_R$ we also require that the event $$\cE_1 := \sbrac{\sum_{i\in[n]} S_i\le K_1n}$$ holds, for a sufficiently large constant~$K_1$. A quick computation verifies that the probability of $\cE_1^\mathrm{c}$ is exponentially small. As discussed in the proof of Proposition \ref{prop:positive}, $\E{\sum_{i\in[n]}S_i}=n\int w\dd\eta_n\le C_2n$, for some constant~$C_2$ independent of~$n$. For~$\theta>0$, independence and~$e^{\theta t}\le1+(e^\theta-1)t$ for~$0\le t\le1$ give
\begin{equation*}
 \E{\exp\!\brac{\theta\sum_{i\in[n]}S_i}}
 \le\exp\!\brac{(e^\theta-1)C_2n}\,.
\end{equation*}
Markov's inequality with~$\theta=1$ and~$K_1>(e-1)C_2+1$ yields~$\Prob{\sum_{i\in[n]}S_i>K_1n}\le O(e^{-n})$. Therefore, without loss of generality, we assume that we are on the event $\cE_0 \cap \cE_1$.

Fix a perfect matching~$P$. For a preimage~$M$ with~$\mathsf{P}(M)=P$, let~$I\subseteq[n]$ be the indices of the rows in which~$M$ does not use the edge of~$P$. Then~$|I|=|P\setminus M|$, and the cancellation outside~$I$ is exactly
\begin{equation*}
 \frac{\weight_\la(M)}{\weight_\la(P)}
 =\exp\!\brac{\la\sum_{e\in P\setminus M}x_e}
   \prod_{e\in M\setminus P}w(x_e)\,.
\end{equation*}
By \lemref{lem:augment}, the reciprocal weights of the missing~$P$-edges $\exp(\la \sum_{e \in P\setminus M}x_e)$ contributes at most
\begin{equation*}
 \exp\!\brac{\la K_0j\log^{1+1/d}(en/j)}.
\end{equation*}
Note that we can bound the factor $\prod_{e \in M\setminus P}w(x_e)$ by overcounting the possible weight contributions since each row $i\in I$  contributes a factor of $1$ if it is unmatched, or $w(x_{i,\ell})$ if it is matching to column $\ell$. Thus, for fixed~$I$, the sum of the weight ratios is at most
\begin{equation}\label{eq:totalcosts}
 \exp\!\brac{\la K_0j\log^{1+1/d}(en/j)}
 \prod_{i\in I}(1+S_i)\,.
\end{equation}

It remains to choose the changed rows. Set~$r=\min\{n,\lfloor K_0 j \log(en/j)\rfloor\}$ and~$\vartheta=r/n$. By Lemma \ref{lem:augment}, we know that $|I|\le r$, so the inequality~$\vartheta^{|I|}\ge\vartheta^r$ gives
\begin{align*}
 \sum_{|I|\le r}\prod_{i\in I}(1+S_i)
 &\le\vartheta^{-r}\prod_{i=1}^n\bigl(1+\vartheta(1+S_i)\bigr)\\
 &\le\exp\bigl(r\log(n/r)+(1+K_1)r\bigr)\,.
\end{align*}
The second inequality uses~$\log(1+t)\le t$ and~$\sum_{i\in[n]}S_i\le K_1n$. Since~$j\le r\le K_0j\log(en/j)$ on the event $\cE_0$, the total contribution of this factor is $O(\exp(j\log^2(en/j)))$, including when~$r=n$. Combining this bound with \eqref{eq:totalcosts} and summing over~$P$ recovers $\cZ_{n,0}$, proving ~\eqref{dim:completion-bound}.
\end{proof}

\ssec{From partial to perfect matchings}\label{sec:removal}

The exponent in~\eqref{dim:completion-bound} is sublinear in~$n$ when~$k=o(n)$. For a positive fraction of unmatched vertices, the factor~$z^{2k}$ can instead be made arbitrarily small. Splitting according to these two ranges is the entire reason the monomers can be removed uniformly in graph size.

\begin{prop}[Removing monomers on the finite graphs]\label{prop:perfect}
The limit~$p(\la):=\lim_{z\downarrow0}p(\la, z)$ is finite and equals the almost sure limit as $n\to\infty$ of~$p_n(\la)$. If~$\E{C^q}<\infty$ for some~$q>0$, the pressures are eventually integrable and the convergence also holds in~$L^1$.
\end{prop}

\begin{proof}
Let~$\cA_n$ be the events from Lemmas~\ref{lem:augment} and~\ref{dim:completion} such that ~$\Prob{\cA_n^\mathrm{c}}=O(n^{-c})$ for some $c > 0$. On~$\cA_n$, the identity~$\cZ_{G_n}(z)=\sum_{k=0}^n z^{2k}\cZ_{n,k}$ and bound 
$$
z^{2k}\frac{\cZ_{n,k}}{\cZ_n(\la)} \leq \exp(Ck\log^{\al_d}(en/k) + 2k\log z) \leq \exp(n\varepsilon(z))
$$
where $\varepsilon(z)$ is defined as
$$
    \varepsilon(z)
 =\sup_{0<u\le1}\{Cu\log^{\al_d}(e/u)+2u\log z\}
$$
gives, for~$0<z\le1$,
\begin{align}\label{eq:sandwich}
 0\le p_n(\la, z)-p_n(\la)
 &= \frac{1}{n}\log(\frac{\cZ_{G_n}(z)}{\cZ_n(\la)}) \\
 &=\frac{1}{n}\log(1 + \sum_{k=1}^n z^{2k}\frac{\cZ_{n,k}}{\cZ_n(\la)}) \\
 &\leq \frac{1}{n}\log((n+1)\exp(n\varepsilon(z))) \\
 &\lesssim\frac{\log(n)}n+\varepsilon(z).\label{eq: final}
\end{align}
Fix~$0<\delta<1$. We choose ~$z$ so small that~$C\log^{\al_d}(e/\delta)+2\log z\le0$. Thus, for every $\delta\leq u \leq 1$, the term $Cu\log^{\al_d}(e/u)+2u\log z$ contributes at most $0$. On~$0\le u\le\delta$, we may discard the nonpositive monomer term to obtain
\begin{equation*}
 0\le\varepsilon(z)
 \le C\sup_{0<u\le\delta}u\log^{\al_d}(e/u)\,.
\end{equation*}
The last bound tends to zero with~$\delta$. Therefore, since $\delta$ was arbitrary, we may assume that $\lim_{z\downarrow 0} \varepsilon(z) = 0$. By Borel--Cantelli, the probabilities $\Prob{\cA_n^{\mathrm{c}}}$ are summable and hence $\cA_n$ holds eventually almost surely. Intersect $\sbrac{\cA_n \text{ e.v.}}$ with the almost sure convergence results for monomer penalties $z \in \mathbb{Q} \cap (0,1]$ corresponding to the event
$$
    \cB = \bigcap_{z \in \mathbb{Q} \cap (0,1]} \sbrac{\lim_{n\to\infty}p_{n}(\la,z) = p(\la,z)}. 
$$
By Proposition \ref{prop:positive}, $\cB$ occurs almost surely. Let $(z_k)_{k\geq 1} = (1/k)_{k\geq1} \downarrow 0$ be a decreasing sequence in $\mathbb{Q}\cap (0,1]$. For all such $z = z_k$, taking the $\liminf_{n\to\infty}$ and $\limsup_{n\to\infty}$ of equation \eqref{eq: final}, we get almost surely
$$
    \lim_{n\to\infty} p_n(\la, z) - \varepsilon(z) \leq \liminf_{n\to\infty} p_n(\la) \leq \limsup_{n\to\infty} p_n(\la) \leq \limsup_{n\to\infty}p_n(\la,z) = p(\la,z),
$$
where we used the monotonicity of $p_n(\la) \leq p_n(\la,z)$ for $z > 0$ in the rightmost inequality. Note that by applying Lemma \ref{lem:augment} to the case of the empty matching when $j = n$, there exists a perfect matching $P$ with total scaled length at most $K_0 n$, which implies that $p_n(\la) \geq -\la K_0$ for large enough $n$. Therefore, $p_n(\la, z) \geq p_n(\la) \geq -\la K_0$ is bounded from below for large enough $n$ and thus so is $p(\la, z)$. Moreover, $p(\la,z)$ is nondecreasing in $z$, and hence $\lim_{z\downarrow 0}p(\la,z)$ exists. Taking $k\to\infty$, we have that almost surely,
$$
    \lim_{n\to\infty} p_n(\la) = \lim_{z\downarrow 0}p(\la,z). 
$$

To prove $L^1$ convergence, assume now~$M_q:=\E{C^q}<\infty$ for some~$q>0$. A uniformly bounded second moment implies uniform integrability by~$|X|\I{\{|X|>T\}}\le X^2/T$. Dropping the column restrictions yields~$\cZ_n(\la)\le\prod_iS_i$, hence
\begin{equation*}
 p_n^+\le\frac1n\sum_iS_i\,.
\end{equation*}
The right side has bounded second moment uniformly in~$n$: its mean is bounded, and independence and~$w^2\le w$ bound its variance by~$n^{-2}\sum_{i,j}\E{w(x_{i,j})}=O(n^{-1})$. This makes the positive parts uniformly integrable.

On~$\cA_n$, the empty-matching argument bounds~$p_n^-$ by~$\la K_0$. To control~$\cA_n^c$ even when~$\E{C}=\infty$, choose an integer~$\el$ with~$q\el>2$. For~$n\ge\el$, take the~$\el$ assignments given by cyclic shifts of column labels by~$0,1,\ldots,\el-1$ modulo~$n$. They use disjoint matrix entries, so their original total costs~$T_1,\ldots,T_\el$ are independent. For~$u>0$, a sum of~$n$ nonnegative terms exceeding~$u$ has one term exceeding~$u/n$. Consequently
\begin{equation*}
 \Prob{T_j>u}\le n\Prob{C>u/n}
 \le M_qn^{1+q}u^{-q}\,.
\end{equation*}
Define~$T_*:=\min_{1\le j\le\el}T_j$ and~$u_0:=M_q^{1/q}n^{1+1/q}$. Independence gives
\begin{equation*}
 \Prob{T_*>u}\le\min\{1,(u_0/u)^{q\el}\}\,.
\end{equation*}
Integrating~$\E{T_*^2}=2\int_0^\infty u\Prob{T_*>u}\dd u$ on~$(0,u_0)$ and~$(u_0,\infty)$ yields
\begin{equation*}
 \|T_*\|_{L^2}\le\brac{\frac{q\el}{q\el-2}}^{1/2}u_0\,.
\end{equation*}
At least one of these assignments contributes~$\exp(-\la(\ka n)^{1/d}T_*)$ to~$\cZ_n(\la)$, so
\begin{equation*}
 p_n^-\le\frac{\la(\ka n)^{1/d}}n T_*\,,
 \qquad \|p_n^-\|_{L^2}=O(n^{1/d+1/q})\,.
\end{equation*}
This already proves integrability for~$n\ge\el$. Choose the completion failure exponent also to satisfy~$c>2/d+2/q$. Cauchy--Schwarz then gives
\begin{equation*}
 \E{p_n^-\I{\cA_n^c}}
 \le\|p_n^-\|_{L^2}\Prob{\cA_n^c}^{1/2}
 =O(n^{1/d+1/q-c/2})\larr0\,.
\end{equation*}
Thus the negative parts are bounded on~$\cA_n$ and have vanishing expected contribution off~$\cA_n$. Together with the positive-part bound, this proves uniform integrability for~$n\ge\el$. Almost-sure convergence and this tail condition imply~$\E{|p_n(\la)-p(\la)|}\to0$.
\end{proof}
\section{The limiting cavity law}\label{sec:hardlaw}

Section \ref{sec:soft} and Section \ref{sec:completion} together nearly complete the final limiting pressure claimed in Theorem \ref{thm:main}. However, the expression still only reaches $\lim_{z\downarrow0}p(\la, z)$, not the originally-claimed limit. Our tasks in Section \ref{sec:hardlaw} are two-fold: to show the existence of a limiting message $\cH$ as $z\downarrow 0$ and that every subsequential limit of $(\cH_{z})_{0< z \leq 1}$ has the same unique law satisfying \eqref{dim:rde}. 
\ssec{Uniform tails for the logarithmic message}\label{dim:messages}
For the first task, we shall use the logarithmic message
\begin{equation*}
 U_{\la,z}:=\la^{-1}\log \cH_z\,
\end{equation*}
and show that its tail bound controls both the location of~$\cH_z$ and rare large child messages in the infinite Poisson sum. 
\begin{lem}[Log-message tails]\label{dim:log-tails}
For every~$d>0$ and~$\la_0>0$, there are constants~$c,K>0$, depending only on~$d,\la_0$, such that
\begin{equation*}
 \Prob{|U_{\la,z}|>t}\le Ke^{-ct^d},
 \qquad t\ge0,\quad \la\ge\la_0,\quad 0<z\le1\,.
\end{equation*}
\end{lem}

\begin{proof}
Write~$U=U_{\la,z}$. On one realization of the recursion define the root output~$\cH_{\varnothing}=(z+\sum_iw(\xi_i)\cH_{z,i})^{-1}$ and its incident-edge probabilities~$a_i=w(\xi_i)\cH_{z,i}\cH_{\varnothing}$. They satisfy~$\sum_i a_i=1-z\cH_{\varnothing}\le1$. Formula~\eqref{eq:mecke} therefore gives
\begin{equation*}
 \int_0^\infty\E{\frac{e^{-\la x}\cH_z\cH_z'}{1+e^{-\la x}\cH_z\cH_z'}}\,\nu_d(\dd x)
 =1-z\E{\cH_z}\le1\,.
\end{equation*}
Here~$\cH_z'$ is an independent copy of~$\cH_z$; put~$U'=\la^{-1}\log \cH_z'$. If~$U>s$,~$U'>s$, and~$x\le2s$, the fraction is at least~$1/2$. Hence
\begin{equation*}
 \tfrac12(2s)^d\Prob{U>s}^2\le1\qquad(s>0)\,.
\end{equation*}
In particular, the probability of~$U>s$ tends uniformly to zero as~$s\to\infty$. This first bound controls upper quantiles; it does not yet give the asserted exponential decay in~$t^d$.

\emph{The median stays bounded below.}
A median~$b$ of~$U$ satisfies both~$\Prob{U\le b}\ge1/2$ and~$\Prob{U\ge b}\ge1/2$. Retain the children with~$U_i:=\la^{-1}\log \cH_{z,i}\ge b$. Independent marking makes their lengths a Poisson process of intensity~$\Prob{U\ge b}\nu_d$. Its first length~$\xi_*$ therefore satisfies
\begin{equation*}
 \Prob{\xi_*>x}\le e^{-x^d/2}\qquad(x\ge0)\,.
\end{equation*}
Keeping this single child in the denominator gives~$\la^{-1}\log \cH_{\varnothing}+b\le\xi_*$. Since~$\la^{-1}\log \cH_{\varnothing}$ has the same law as~$U$,
\begin{equation}\label{eq:median-tail}
 \Prob{U+b>t}\le e^{-t^d/2}\qquad(t\ge0)\,.
\end{equation}

Suppose that a sequence of admissible parameters has~$b\to-\infty$. Upper tightness implies~$U+b\to-\infty$ in probability;~$\la\ge\la_0$ then gives~$e^{\la(U+b)}\to0$ in probability. Rewrite the recursion, with~$U$ now the logarithmic root output, as
\begin{equation*}
 e^{\la(b-U)}=ze^{\la b}
   +\sum_i e^{-\la\xi_i}e^{\la(U_i+b)}\,.
\end{equation*}
The first term tends to zero. At a fixed cutoff~$R$, the sum with~$\xi_i\le R$ also tends to zero in probability: its number of terms has the fixed Poisson mean~$R^d$, and each marked factor tends to zero in probability.

For the remaining sum, the probability of any point with~$U_i+b>\xi_i/2$ and~$\xi_i>R$ is at most~$\int_R^\infty e^{-(x/2)^d/2}\nu_d(\dd x)$ by~\eqref{eq:median-tail}. On the complementary event each summand is at most~$e^{-\la_0\xi_i/2}$. For any~$\epsilon>0$, Markov's inequality bounds the probability that this remaining sum exceeds~$\epsilon$ by
\begin{equation*}
 \int_R^\infty e^{-(x/2)^d/2}\nu_d(\dd x)
 +\epsilon^{-1}\int_R^\infty e^{-\la_0x/2}\nu_d(\dd x)\,.
\end{equation*}
This tends to zero uniformly along the parameter sequence as~$R\to\infty$. Taking the parameter limit first at fixed~$R$, then increasing~$R$, makes the entire right side of the rewritten recursion tend to zero in probability. But its left side is at least one on~$\{U\le b\}$, an event of probability at least~$1/2$. This contradiction proves a uniform lower bound~$b\ge-B_0$, where~$B_0$ depends only on~$d,\la_0$.

For~$t\ge2B_0$,~\eqref{eq:median-tail} now gives~$\Prob{U>t}\le\exp(-t^d/2^{d+1})$. Increasing its prefactor to cover smaller~$t$ proves~$\Prob{U>t}\le Ke^{-ct^d}$.

\emph{The lower tail.}
The expected number of children with~$U_i>(t+\xi_i)/2$ is at most
\begin{equation*}
 \int_0^\infty Ke^{-c((t+x)/2)^d}\nu_d(\dd x)
 \le K'e^{-c't^d}\,.
\end{equation*}
The last estimate uses~$(t+x)^d\ge\min(1,2^{d-1})(t^d+x^d)$, which holds for every~$d>0$. If none of these exceptional children occurs, then
\begin{equation*}
 \cH_{\varnothing}^{-1}\le1+e^{\la t/2}S_{\la_0},
 \qquad S_{\la_0}:=\sum_i e^{-\la_0\xi_i/2}\,.
\end{equation*}
This last variable has every positive exponential moment. Indeed, the Poisson exponential formula~\cite[Theorem~3.9 and Exercise~3.4]{LastPenrose} gives, for~$\theta>0$,
\begin{equation*}
 \E{e^{\theta S_{\la_0}}}
 =\exp\!\brac{\int_0^\infty
   (e^{\theta e^{-\la_0x/2}}-1)\,\nu_d(\dd x)}<\infty\,.
\end{equation*}
Finiteness follows from~$e^{\theta y}-1\le(e^\theta-1)y$ for~$0\le y\le1$ and the integrability of~$e^{-\la_0x/2}$ against~$\nu_d$.

For sufficiently large~$t$, the event~$\la^{-1}\log \cH_{\varnothing}<-t$ on the preceding nonexceptional event forces~$S_{\la_0}>e^{\la_0t/2}/2$. Markov's inequality applied to~$e^{S_{\la_0}}$ gives
\begin{equation*}
 \Prob{S_{\la_0}>e^{\la_0t/2}/2}
 \le\E{e^{S_{\la_0}}}\exp(-e^{\la_0t/2}/2)
 \le Ke^{-ct^d}
\end{equation*}
after adjusting the constants. Add the exceptional-child probability to obtain the lower-tail estimate.
\end{proof}

\ssec{Existence of $\cH$ when $z = 0$}
\begin{prop}[Existence at $z=0$]\label{prop:existence z = 0} Fix $d > 0$ and let $(\xi_i)_{i\geq 1}$ be a Poisson point process on $(0,\infty)$ with intensity $\nu_d(\mathrm{d}x) = dx^{d-1}\mathrm{d}x$. Let $\cH_{\la,z}$ be the positive-monomer message as defined in \eqref{eq:positive} i.e.
$$
    \cH_{\la,z} \overset{\mathrm{d}}{=} \brac{z + \sum_{i\geq 1}e^{-\la \xi_i}\cH_{\la,z,i}}^{-1},
$$
where $\cH_{\la,z,i}$ are independent copies of $\cH_{\la,z}$, independent of $(\xi_i)_{i\geq 1}$. If $z_k \downarrow 0$ and $\la_k\to \la > 0$, then the laws 
of $(\cH_{\la_k,z_k})_{k\geq1}$ are relatively compact in the weak topology on $\mathscr{P}((0,\infty))$. In particular, every subsequential limit $\cH$ of $(\cH_{\la_k,z_k})_{k\geq 1}$ satisfies 
$$
    \cH \overset{\mathrm{d}}{=} \brac{\sum_{i\geq 1} e^{-\la \xi_i}\cH_{i}}^{-1},
$$
where $\cH_{i}$ are independent copies of $\cH$, independent of $(\xi_{i})_{i\geq 1}$ and the sum is both positive and finite almost surely. 
\end{prop}
\begin{proof}
Discard finitely many terms so that $0<z_k\le1$ and $0<\la_-\le\la_k\le\la_+<\infty$. Set $H_k=\cH_{\la_k,z_k}$ and $U_k=\la_k^{-1}\log H_k$. By~\lemref{dim:log-tails}, $\sup_k\Prob{|U_k|>t}\le Ke^{-ct^d}$ for $t\ge0$. Since $\la_k$ is bounded, $\log H_k=\la_k U_k$ is tight, and hence so is $H_k$ on $(0,\infty)$. Along any subsequence with $H_k\Rightarrow\cH$, we have $U_k\Rightarrow U=\la^{-1}\log\cH$ because $\la_k\to\la$. The same log-tail bound holds for $U$ by Portmanteau.

Let $S_k=\sum_i e^{-\la_k\xi_i}H_{k,i}$ and $S=\sum_i e^{-\la\xi_i}\cH_i$, with independent marks of the respective laws. Write $S_{k,R}$ and $S_R$ for their truncations to $\xi_i\le R$. Conditional on the finitely many retained lengths, weak convergence of the marks and $\la_k\to\la$ give $S_{k,R}\Rightarrow S_R$. Let $N_{k,R}$ count points with $\xi_i>R$ and $\la_k^{-1}\log H_{k,i}>\xi_i/2$. By Campbell's formula,
\begin{equation}\label{eq: N_kr}
    \Prob{N_{k,R}\geq 1} \leq \E{N_{k,R}} = \int_R^{\infty} \Prob{U_k > x/2}\nu_d(\mathrm{d}x) \leq \int_{R}^\infty K e^{-c(x/2)^d}\nu_d(\mathrm{d}x). 
\end{equation}
On the event $\sbrac{N_{k,R}=0}$, every point with $\xi_i>R$ satisfies $e^{-\la _k \xi_i}H_{k,i}\leq e^{-\la_- \xi_i/2}$ and therefore,
$$
    S_k - S_{k,R} \leq \sum_{i\geq 1: \xi_i > R} e^{-\la_-\xi_i/2}. 
$$
By~\eqref{eq: N_kr} and Markov's inequality, for every $\epsilon>0$,
\begin{equation}\label{eq:poisson-tail}
 \sup_k\Prob{S_k-S_{k,R}>\epsilon}
 \le \int_R^\infty Ke^{-c(x/2)^d}\nu_d(\dd x)
   + \epsilon^{-1}\int_R^\infty e^{-\la_- x/2}\nu_d(\dd x)
 \xrightarrow[R\to\infty]{}0.
\end{equation}
For the limiting marks, the tail bound for $U$ also makes the expected number of points with $\la^{-1}\log\cH_i>\xi_i/2$ finite. The other terms of $S$ are bounded by the summable $e^{-\la\xi_i/2}$, so $S<\infty$ almost surely. Also $S>0$ since the process has infinitely many points and all marks are positive. Now $S_{k,R}\Rightarrow S_R$ at fixed $R$, $S_R\uparrow S$, and the uniform tail estimate yields $S_k\Rightarrow S$. Therefore $(z_k+S_k)^{-1}\Rightarrow S^{-1}$ by continuity at positive $S$. The law of $(z_k+S_k)^{-1}$ is the law of $H_k$ by~\eqref{eq:positive}, so $\cH\stackrel{\mathrm d}=S^{-1}$.
\end{proof}
Therefore, Proposition \ref{prop:existence z = 0} states that for every sequence~$z_k\downarrow0$, $(\cH_{z_k})_{k\geq 1}$ has a weakly converging subsequence such that (up to subsequences) $\cH_{z_k}\Rightarrow \cH$ for a positive finite random variable~$\cH$ satisfying the recursive distributional equation \eqref{dim:rde}. However, existence alone does not identify the whole family~$(\cH_z)_{z \in (0,1]}$. In this next section, we show that any two positive finite solutions of~\eqref{dim:rde} must agree.

% We must pass through the infinite Poisson sum, not merely a finite recursion. For a cutoff~$R$, define its omitted part by
% \begin{equation*}
%  T_{z,R}=\sum_{\xi_i>R}e^{-\la\xi_i}\cH_{z,i}\,.
% \end{equation*}
% For each child put~$U_{\la,z,i}=\la^{-1}\log \cH_{z,i}$. Split according to whether some~$U_{\la,z,i}>\xi_i/2$. \lemref{dim:log-tails} and Markov's inequality give, for~$\epsilon>0$,
% \begin{equation}\label{eq:poisson-tail}
%  \sup_{0<z\le1}\Prob{T_{z,R}>\epsilon}
%  \le\int_R^\infty Ke^{-c(x/2)^d}\nu_d(\dd x)
%  +\epsilon^{-1}\int_R^\infty e^{-\la x/2}\nu_d(\dd x)
%  \larr0\,.
% \end{equation}
% For fixed~$R$, the retained sum has finitely many terms, so weak convergence of its independent marks gives convergence of the sum in distribution. The estimate above removes~$R$ uniformly. The limiting full sum with marks~$\cH_i$ is finite: there are almost surely only finitely many exceptional marks~$\la^{-1}\log \cH_i>\xi_i/2$, and the other terms are bounded by the summable~$e^{-\la\xi_i/2}$. It is strictly positive because the process has infinitely many points and the marks are positive.

% Only after passing to this positive full sum do we take its reciprocal. A finite-cutoff sum can equal zero, so inversion at that earlier stage would not be justified. The recursion~\eqref{eq:positive} now yields~\eqref{dim:rde} for the subsequential limit~$\cH$. This proves existence. The same reasoning permits a convergent sequence of positive temperatures~$\la_j\to\la$ by using one common positive lower bound~$\la_0$.

\ssec{Comparing the edge choices}

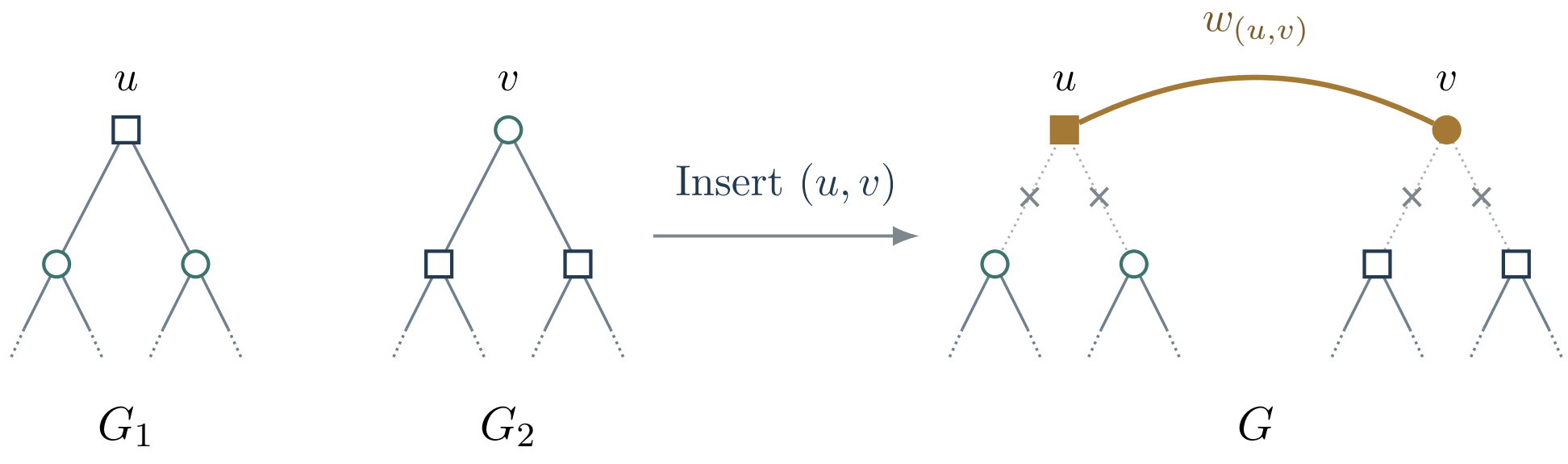
\begin{figure}[!htbp]
\centering
\begin{tikzpicture}[x=.95cm,y=1cm,font=\small,
 row/.style={rectangle,draw=figNavy,fill=white,line width=.8pt,minimum size=6pt,inner sep=0pt},
 col/.style={circle,draw=figTeal,fill=white,line width=.8pt,minimum size=6pt,inner sep=0pt},
 branch/.style={draw=figNavy!65,line width=.65pt},
 blocked/.style={draw=figGray!65,densely dotted,line width=.55pt},
 cross/.style={draw=figGray,line width=.75pt},
 occupied/.style={draw=figOchre,fill=figOchre},
 graph name/.style={font=\normalsize,text=black,anchor=north}]
 % The lower branches are identical on both sides.
 \foreach \x in {1,4.3,9.1,12.4}{
   \foreach \dx in {-.6,.6}{
     \draw[branch] (\x+\dx-.28,.62) -- (\x+\dx,1.15) -- (\x+\dx+.28,.62);
     \draw[branch,densely dotted] (\x+\dx-.28,.62) -- (\x+\dx-.4,.37);
     \draw[branch,densely dotted] (\x+\dx+.28,.62) -- (\x+\dx+.4,.37);
   }
 }
 \foreach \x in {1,4.3}{
   \draw[branch] (\x-.6,1.15) -- (\x,2.25) -- (\x+.6,1.15);
 }
 % Occupying the inserted edge blocks the other edges at its endpoints.
 \foreach \x in {9.1,12.4}{
   \foreach \dx in {-.6,.6}{
     \draw[blocked] (\x,2.25) -- (\x+\dx,1.15);
     \draw[cross] (\x+.5*\dx-.075,1.625) -- (\x+.5*\dx+.075,1.775);
     \draw[cross] (\x+.5*\dx-.075,1.775) -- (\x+.5*\dx+.075,1.625);
   }
 }
 \foreach \x in {.4,1.6,8.5,9.7}{\node[col] at (\x,1.15) {};}
 \foreach \x in {3.7,4.9,11.8,13}{\node[row] at (\x,1.15) {};}
 \node[row,label={[label distance=2pt]above:$u$}] at (1,2.25) {};
 \node[col,label={[label distance=2pt]above:$v$}] at (4.3,2.25) {};
 \node[row,occupied,label={[label distance=2pt]above:$u$}] (uAfter) at (9.1,2.25) {};
 \node[col,occupied,label={[label distance=2pt]above:$v$}] (vAfter) at (12.4,2.25) {};
 \draw[figOchre,line width=1.25pt]
   (uAfter) to[bend left=26] node[midway,above=3pt,text=figOchre!72!black] {$w_{(u,v)}$} (vAfter);
 \draw[-{Latex[length=2.2mm]},figGray,line width=.75pt] (5.55,1.38) -- (7.85,1.38)
   node[midway,above=3pt,text=figNavy] {Insert $(u,v)$};
 \node[graph name] at (1,.11) {$G_1$};
 \node[graph name] at (4.3,.11) {$G_2$};
 \node[graph name] at (10.75,.11) {$G$};
\end{tikzpicture}
\caption{Inserting $(u,v)$ joins the rooted graphs $G_1$ and $G_2$ into $G$. The inserted edge is shown occupied, so the crossed branch edges at $u$ and $v$ are unavailable to the matching. Matchings omitting or using $(u,v)$ have relative weights $1$ and $w_{(u,v)}\cA\cB$, respectively, after normalizing by the product of the separate partition functions. The branches are schematic.}\label{fig:insertion}
\end{figure}

For $z>0$, consider disjoint finite graphs $G_1$ and $G_2$ with rooted vertices $u$ and $v$ respectively. Set $\cA = {\cZ_{G_1 - u}(z)}/{\cZ_{G_1}(z)}$ and $\cB = {\cZ_{G_2 - v}(z)}/{\cZ_{G_2}(z)}$. Let $G$ be the graph with edge $e=(u,v)$ inserted with weight $w_{(u,v)}$ (see Figure \ref{fig:insertion}). Notice that $\cZ_G(z)$ can be written as
$$
    \cZ_G(z) = \cZ_{G_1}(z) \cZ_{G_2}(z)(1 + w_{(u,v)}\cA \cB)
$$
and hence the probability of $(u,v)$ belonging to the matching under the Gibbs measure $\mu_{G,z}$ is $\frac{w_{(u,v)}\cA\cB}{1 + w_{(u,v)}\cA \cB}$. Let $\cA^{+e} = \cH_{G,z}(u) = \cA/(1+w_{(u,v)}\cA\cB)$ and $\cB^{+e} = \cH_{G,z}(v) = \cB/(1 + w_{(u,v)}\cA\cB)$ be the \textit{deletion} \textit{messages} at $u$ and $v$, respectively, in the joined graph $G$. We will use the corresponding tree identities to identify the limiting cavity law for the overlap.

At $z=0$, take $\cA,\cB$ to be the cavity messages at the endpoints of a tree edge~$e$ of weight~$w$, with that edge omitted. Restoring the edge adds~$w\cB$ to~$\cA^{-1}$ and~$w\cA$ to~$\cB^{-1}$ in the two tree recursions. The endpoint messages with the edge present are therefore~$(\cA^{-1}+w\cB)^{-1}$ and~$(\cB^{-1}+w\cA)^{-1}$. Write~$\cA^{+e},\cB^{+e}$ for these messages and set~$a=w\cB\cA^{+e}$, the coefficient of this edge in the recursion at the first endpoint. Then
\begin{equation}\label{eq:insertion}
 a=\frac{w\cA\cB}{1+w\cA\cB},\qquad
 \cA^{+e}=(1-a)\cA,\qquad
 \cB^{+e}=(1-a)\cB\,.
\end{equation}

The insertion identity~\eqref{eq:insertion} gives~$a(1-a)=w\cA^{+e}\cB^{+e}$. Taking logarithms separates the two endpoints. When two candidate laws are compared on the same edge, its weight cancels, leaving one log-message difference at each endpoint. This explains the comparison quantity in the following lemma.

\begin{lem}[A probability-vector inequality]\label{lem:simplex}
Let~$a=(a_i)$ and~$\wt a=(\wt a_i)$ be finite or countable sequences with~$a_i>0$,~$\wt a_i>0$, and~$\sum_i a_i=\sum_i\wt a_i=1$. Assume that their common index set has at least three elements. Then
\begin{equation}\label{eq:simplex}
 \cD(a,\wt a):=
 \sum_i(a_i-\wt a_i)
 \log\frac{a_i(1-a_i)}{\wt a_i(1-\wt a_i)}\ge0\,,
\end{equation}
with equality if and only if~$a=\wt a$. For countably many coordinates, the sum is allowed to be~$+\infty$.
\end{lem}

The proof is given in Appendix~\ref{app:simplex}.

\begin{lem}[An invariant coupling]\label{lem:invariant-coupling}
Fix $d,\la>0$, and let $\mu$ and $\tilde\mu$ be two solution laws of
\eqref{dim:rde}. There is a coupling $\rho\in\Pi(\mu,\tilde\mu)$ such that
if $(\cH_i,\wt\cH_i)_{i\geq1}$ are $\iid$ pairs of law~$\rho$, independent
of the Poisson points $(\xi_i)_{i\geq1}$ in \eqref{dim:rde}, then
\[
 \rho = \cL\brac{\brac{\sum_{i\geq1}e^{-\la\xi_i}\cH_i}^{-1},
       \brac{\sum_{i\geq1}e^{-\la\xi_i}\wt\cH_i}^{-1}}.
\]
\end{lem}

\begin{proof}
Let $(\mu,\tilde\mu) \in \mathscr{P}((0,\infty))^2$ be two probability measures satisfying \eqref{dim:rde}. Given $\rho\in\Pi(\mu,\tilde\mu)$, mark the same Poisson points by independent pairs $(\cH_i,\wt\cH_i)$ of law~$\rho$. Let $\mathfrak T(\rho)$ be the joint law of the two root outputs
\[
 \mathfrak{T}(\rho) = \cL\brac{\brac{\sum_{i\geq 1} w(\xi_i)\cH_i}^{-1}, \brac{\sum_{i\geq 1} w(\xi_i)\wt\cH_i}^{-1}}.
\]
By \eqref{dim:rde}, it is clear that $\mathfrak T: \Pi(\mu,\tilde\mu) \to \Pi(\mu,\tilde{\mu})$. Moreover, the set $\Pi(\mu,\tilde\mu)$ is nonempty, convex, and weakly compact by tightness of the fixed marginals and closedness of the marginal constraints \cite[Lemma 4.4]{villani2009optimal}. We show that $\mathfrak{T}$ is continuous on the space of couplings equipped with the weak topology. For any $\rho_m \Rightarrow \rho$, let $k > 0$ be a fixed integer and $\{(\cH_{m,i}, \tilde{\cH}_{m,i})\}_{i\geq 1}$ be $\iid$ draws from $\rho_m$. Define the partial sum couplings
$$\mathbf{S}_m := \brac{\sum_{i\geq 1}w(\xi_i)\cH_{m,i}, \sum_{i\geq 1}w(\xi_i)\tilde{\cH}_{m,i}}, \quad \mathbf{S}_{m}^{(k)} = \brac{\sum_{i=1}^{k}w(\xi_i)\cH_{m,i}, \sum_{i=1}^{k} w(\xi_i)\tilde{\cH}_{m,i}}.$$
For fixed $k$, clearly $\mathbf{S}_m^{(k)}\Rightarrow \mathbf{S}^{(k)},$ where $\mathbf{S}^{(k)}$ is the corresponding $k$-level sum for $\rho$. Since the laws of the discarded tail $\mathbf{S}_m - \mathbf{S}_m^{(k)}$ is identical for all $m$, and the infinite sum is finite almost surely by Proposition \ref{prop:existence z = 0}, for any $\epsilon > 0$, as $k \to \infty$,
$$
    \sup_{m \geq 1} \Prob{\norm{\mathbf{S}_m - \mathbf{S}_m^{(k)}}_{\infty}\geq \epsilon} \longrightarrow 0.
$$
The same holds for the limiting coupling $\rho$ as well. Therefore, for any $1$-Lipschitz test function $f$ with $\norm{f}_\infty \leq 1$,
$$
\abs{\E{f(\mathbf{S}_m)} - \E{f(\mathbf{S})}} \leq \abs{\E{f(\mathbf{S}_{m}^{(k)})} - \E{f(\mathbf{S}^{(k)})}} + 2\epsilon + o_{\epsilon}(1).
$$
Taking $m\to\infty$ for fixed $k$, $k\to \infty$, and $\epsilon \downarrow 0$, we arrive at $\mathbf{S}_m \Rightarrow \mathbf{S}$. Again, by Proposition \ref{prop:existence z = 0}, the coordinates of $\mathbf{S}$ are positive almost surely, so by the continuous mapping theorem $\mathfrak{T}(\rho_m) \Rightarrow \mathfrak{T}(\rho)$. Since the weak topology on finite signed measures, induced by bounded continuous test functions, is Hausdorff and locally convex, by the Schauder--Tychonoff fixed point theorem, there exists a coupling $\rho \in \Pi(\mu, \tilde{\mu})$ such that $\rho=\mathfrak T(\rho)$~\cite[Theorem~3.2]{Bonsall}.
\end{proof}

\ssec{Uniqueness by exchanging endpoints}

\begin{prop}[The zero-fugacity law]\label{prop:unique}
For each~$d>0$ and~$\la>0$, there is exactly one probability law on~$(0,\infty)$ solving~\eqref{dim:rde}. If~$\cH$ has this law, then the positive-monomer message~$\cH_{\la,z}$ from~\eqref{eq:positive} converges in distribution to~$\cH$ as~$z\downarrow0$. Moreover, set~$U_{\la,0}=\la^{-1}\log\cH$. For each~$\la_0>0$, the constants~$c,K$ of~\lemref{dim:log-tails} satisfy
\begin{equation*}
 \Prob{|U_{\la,0}|>t}\le Ke^{-ct^d},
 \qquad t\ge0,\quad \la\ge\la_0.
\end{equation*}
\end{prop}

\begin{proof}
Fix $\la>0$. By~\propref{prop:existence z = 0}, at least one solution exists, and every subsequential limit as $z\downarrow0$ is a solution. To compare arbitrary solutions, first note that each has an integrable logarithm. Indeed, for any solution take a median~$b$ of $U=\la^{-1}\log\cH$. The thinning argument of~\eqref{eq:median-tail} with this median and $z=0$ gives $\Prob{U+b>t}\le e^{-t^d/2}$ for $t\ge0$. The lower-tail argument of~\lemref{dim:log-tails} then applies with the resulting upper-tail bound.

Let $\mu$ and $\tilde\mu$ be two solution laws, and take an invariant coupling
$\rho$ given by \lemref{lem:invariant-coupling}. Let $(\cH_i,\wt\cH_i)_{i\geq1}$
be the independent child pairs of law~$\rho$ in the lemma, and write
$(\cH_{\varnothing},\wt\cH_{\varnothing})$ for the corresponding root outputs.

Use this coupling and write $a_i=w(\xi_i)\cH_i\cH_{\varnothing}$, with $\wt a_i$ defined in the second system. Both are probability vectors with infinitely many positive entries. By~\eqref{eq:insertion}, the full child values are $\cH_{v_i}=(1-a_i)\cH_i$ and $\wt\cH_{v_i}=(1-\wt a_i)\wt\cH_i$. Write $\De_v=\log(\cH_v/\wt\cH_v)$ for the full log difference at either endpoint and set $X_i=a_i-\wt a_i$. Then
\[
 \log\frac{a_i(1-a_i)}{\wt a_i(1-\wt a_i)}
   = \De_{\varnothing}+\De_{v_i},
 \qquad \sum_i X_i=0.
\]

Under Mecke's formula~\eqref{eq:mecke}, the cavity pairs at the two ends of an inserted edge are independent with law~$\rho$. One is the root output of the remaining Poisson process, and the other is a fresh mark. Swapping them preserves the edge probabilities and exchanges the full endpoint log differences by~\eqref{eq:insertion}. Hence for every nonnegative measurable~$f$,
\[
 \E{\sum_i f(X_i,\De_{v_i})}=\E{\sum_i f(X_i,\De_{\varnothing})}.
\]
With $f(u,v)=|uv|$ the common expectation is at most $2\E{|\De_{\varnothing}|}<\infty$. Applying the same identity to the signed product and using $\sum_i X_i=0$ gives
\[
 \E{\cD(a,\wt a)}
 =\E{\sum_i X_i\De_{v_i}}
 =\E{\De_{\varnothing}\sum_i X_i}=0.
\]
By~\lemref{lem:simplex}, $a_i=\wt a_i$ for every~$i$ almost surely. Thus each child log ratio $\log(\cH_i/\wt\cH_i)$ equals $-\De_{\varnothing}$. Two independent child log ratios agree almost surely, so their common law is a point mass at some $c$. Invariance gives $\De_{\varnothing}=c$ almost surely, and hence $c=-c=0$. The coupling is supported on equal pairs, proving $\mu=\tilde\mu$.

Relative compactness from~\propref{prop:existence z = 0} and uniqueness give $\cH_{\la,z}\Rightarrow\cH$ as $z\downarrow0$. The bound of~\lemref{dim:log-tails} passes to the limit by Portmanteau on the open log-tail events.
\end{proof}

\ssec{The pressure formula}\label{sec:assembly}
We now present the promised proof of the limiting free energy formula, as well as the remaining arguments needed to prove the remaining Corollaries \ref{cor:energy} and \ref{cor:cost-budget}. 

\begin{proof}[Proof of \thmref{dim:main}]
By \propref{prop:perfect}, almost surely,
$$
    \lim_{n\to\infty}p_n(\la) = p(\la) = \lim_{z\downarrow 0}p(\la,z),
$$
with convergence also in $L^1$ if $\E{C^q} < \infty$ for some $q > 0$. Moreover, by \propref{prop:unique}, there exists a unique random variable $\cH$ in $(0,\infty)$ such that $\lim_{z\downarrow 0}\cH_z = \cH$ almost surely and satisfying \eqref{dim:rde}. As before, setting $U_{\la,z} = \la^{-1}\log \cH_z$, define $S_z := U_{\la,z} + U'_{\la,z}$ where $(U_{\la,z},U'_{\la,z})$ are independent copies. Define 
$$
    F(s) = \int_0^\infty \log(1 + e^{\la(s - x)})\nu_d(\mathrm{d}x).
$$
By Tonelli's theorem, 
$$
    p(\la,z) = -(\E{\la S_z + F(S_z)}). 
$$
Since $S_z \Rightarrow S_0 :=\la^{-1}(\log \cH + \log \cH^\prime)$ for $(\cH, \cH^\prime)$ independent, we need to show $F(S_z) \Rightarrow F(S_0)$. Note that $F(s)$ is continuous. Indeed, for any sequence $(s_k)_k$ converging to $s$, there exists some $M < \infty$ such that $|s_k|\leq M$ and so
$$
    0 \leq d x^{d-1}\log(1+e^{\la(s_k - x)})\leq de^{\la M}e^{-\la x}x^{d-1},
$$
which is in $L^1([0,\infty))$. Therefore, by dominated convergence, $\lim_{k\to\infty}F(s_k) = F(s)$. By another application of the continuous mapping theorem, we have that for any constant $M > 0$, $\E{F(S_z) \wedge M} \to \E{F(S_0) \wedge M}$. Setting $s_+ = \max(s,0)$ and using the inequality $\log(1+e^{x})\leq x_+ + \log 2$, splitting the above integral around $2s_+$ gives
$$
    0\leq F(s) \leq (2s_+)^d(\log 2 +  \la s_+) + \int_0^\infty e^{-\la x/2}\nu_d(\mathrm{d}x) = O(1+|s|^{d+1}).
$$
Therefore, for some universal constant $C > 0$,
$$
     \E{F(S_z)^2} \leq 2C^2\brac{1+\E{|S_z|^{2(d+1)}}}.
$$
Using the log-message tail bound from Lemma \ref{dim:log-tails}, uniformly in $0 < z \leq 1$,
$$
    \E{|S_z|^{2(d+1)}} = \int_0^\infty 2(d+1)t^{2d+1}\Prob{|S_z|>t}\,\mathrm{d}t \leq 4K(d+1)\int_0^\infty t^{2d+1}e^{-c(t/2)^d}\,\mathrm{d}t < \infty.
$$
Consequently, by a standard truncation argument, $\E{F(S_z)\wedge M}\to \E{F(S_z)}$ as $M\to\infty$, uniformly in $z \in [0,1]$, and so 
$$
    \limsup_{z\downarrow 0}\abs{\E{F(S_z)}-\E{F(S_0)}}=0.
$$
A similar argument also proves convergence of expectations for $\E{\la S_z}$ to $\E{\la S_0}$. Therefore,
$$
    \lim_{z\downarrow 0}p(\la,z) = p(\la) = -(\E{\la S_0 + F(S_0)}) = -2\E{\log \cH} - \int_0^\infty \E{\log(1 + e^{-\la x}\cH \cH^\prime)}\nu_{d}(\mathrm{d}x),
$$
as required.
\end{proof}

\begin{lem}[Pressure derivative]\label{lem:pressure-derivative}
The limiting pressure $p(\la)$ is continuously differentiable on $(0,\infty)$. Moreover, $u(\la) :=-p^\prime(\la)$ is given by the cavity integral in \corref{cor:energy} 
$$
    u(\la) = \int_0^\infty x\E{\frac{e^{-\la x}\cH\cH^\prime}{1+e^{-\la x}\cH\cH^\prime}}\nu_d(\mathrm{d}x).
$$
and satisfies $u(\la)\ge d/\la$.
\end{lem}

\begin{proof}
For $\la>0$ and $x,h,k\ge0$, set
\[
 f_\la(x,h,k)=\frac{e^{-\la x}hk}{1+e^{-\la x}hk}.
\]
For $z>0$, write $\cH_{\la,z}$ for the cavity message of
\lemref{lem:positive} and $\cH_{\la,z,R}$ for its counterpart obtained by
restricting the Poisson recursion to lengths at most $R>0$. At $z=0$,
write $\cH_{\la,0}$ for the message of \propref{prop:unique}. Define
\begin{align*}
 Q_{\la,z}(x)
 &=\E{f_\la(x,\cH_{\la,z},\cH'_{\la,z})}, && z\ge0,\\
 Q_{\la,z,R}(x)
 &=\E{f_\la(x,\cH_{\la,z,R},\cH'_{\la,z,R})}, && z>0.
\end{align*}
Let $p_{n,R}(\la,z)=n^{-1}\log\cZ_{G_{n,R}}(z)$, and let
$\mu_{n,R,\la,z}$ be the associated Gibbs measure. Fix $z,R>0$. Differentiating gives
\begin{equation}\label{eq:cutoff-pressure-derivative}
 -\partial_\la p_{n,R}(\la,z)
 =\frac1n\sum_{e\in E(G_{n,R})}x_e\,\mu_{n,R,\la,z}(e\in M).
\end{equation}
The right-hand side is at most $R$, since each matching contains at most
$n$ edges and each retained edge has length at most $R$. This bound also
justifies differentiating under the expectation over the costs.

Fix one possible edge $e=(v,w)$ and put $G:=G_{n,R}\setminus\{e\}$, retaining both vertices $v$ and $w$. Define
\[
 A_n=\frac{\cZ_{G-v}(z)}{\cZ_G(z)},\qquad
 B_n=\frac{\cZ_{G-\{v,w\}}(z)}{\cZ_{G-v}(z)}.
\]
The pair $(A_n,B_n)$ depends only on the lengths of edges other than $e$.
When $e$ is inserted at length $x_e\le R$, the total weights of matchings
omitting and using $e$ are respectively $\cZ_G(z)$ and
$e^{-\la x_e}\cZ_{G-\{v,w\}}(z)$. Thus the occupation probability of the
inserted edge is $f_\la(x_e,A_n,B_n)$. Let $\eta_{n,R}$ be the restriction
of $n\mathcal L(x_e)$ to $[0,R]$.
Exchangeability and independence of $x_e$ from $(A_n,B_n)$ give
\[
 -\partial_\la\E{p_{n,R}(\la,z)}
 =\int_{[0,R]} x\,\E{f_\la(x,A_n,B_n)}\,\eta_{n,R}(\mathrm dx).
\]
We claim that
\[
 (A_n,B_n)\Rightarrow(\cH_{\la,z,R},\cH'_{\la,z,R}),
\]
whose proof we defer.
% At each fixed depth, the two-root exploration of \propref{prop:positive}
% gives joint convergence of the neighborhoods of $v$ and $w$ in $G$ to
% independent cutoff Poisson trees. The two finite neighborhoods are disjoint
% trees with probability tending to one. On this event the neighborhood of
% $w$ is unchanged by deleting $v$. Thus the depth bounds of
% \lemref{lem:positive} apply to $A_n$ on the first neighborhood and to $B_n$
% on the second. On the complementary event we may use $0$ and $1/z$ as
% bounds for both ratios. The two pairs of bounds converge jointly at fixed
% depth. By \lemref{lem:positive} and bounded convergence, the expected
% widths of the two bounding intervals on the limiting trees tend to zero
% as the depth increases. For any bounded Lipschitz test function of the
% pair, the error from replacing the ratios by their lower bounds is at
% most the Lipschitz constant times the sum of the widths. Letting $n$ tend
% to infinity at fixed depth and then sending the depth to infinity proves
% the claim.
The cost assumption \eqref{dim:F} gives $\eta_{n,R}([0,t])\to t^d$ for
$0\le t\le R$. Hence $\eta_{n,R}$ converges weakly as a finite measure to
$\nu_d$ restricted to $[0,R]$. Since $(x,h,k)\mapsto x f_\la(x,h,k)$ is
bounded and continuous on $[0,R]\times[0,1/z]^2$, as $n\to\infty$,
\[
 -\partial_\la\E{p_{n,R}(\la,z)}
 \longrightarrow \int_0^R x\,Q_{\la,z,R}(x)\,\nu_d(\mathrm dx).
\]

Fix $0<a<b<\infty$. Write
$p_R(\la,z)=\lim_{n\to\infty}\E{p_{n,R}(\la,z)}$, as in
\propref{prop:positive}. Integrating the preceding convergence over $[a,b]$
and using the bound $R$ from \eqref{eq:cutoff-pressure-derivative} gives
\[
 p_R(b,z)-p_R(a,z)
 =-\int_a^b\int_0^R x\,Q_{\la,z,R}(x)\,\nu_d(\mathrm dx)\,\mathrm d\la.
\]
At fixed $z>0$, \lemref{lem:positive} gives
$\cH_{\la,z,R}\Rightarrow\cH_{\la,z}$ as $R\to\infty$. Since the messages
are bounded by $1/z$, we have $Q_{\la,z,R}(x)\le z^{-2}e^{-ax}$ for
$\la\in[a,b]$, which is integrable against $x\,\nu_d(\mathrm dx)$.
Using $p_R(\la,z)\to p(\la,z)$ from \propref{prop:positive} and applying
dominated convergence therefore gives
\begin{equation}\label{eq:soft-pressure-increment}
 p(b,z)-p(a,z)
 =-\int_a^b\int_0^\infty x\,Q_{\la,z}(x)\,\nu_d(\mathrm dx)\,\mathrm d\la.
\end{equation}

To pass to $z=0$, write $U_{\la,z}=\la^{-1}\log\cH_{\la,z}$ and $U'_{\la,z}=\la^{-1}\log\cH'_{\la,z}$ for $z\ge0$.
On the event $U_{\la,z},U'_{\la,z}\le x/4$, we have
\[
 e^{-\la x}\cH_{\la,z}\cH'_{\la,z}
 =e^{\la(U_{\la,z}+U'_{\la,z}-x)}\le e^{-ax/2}.
\]
On the complementary event, the kernel is bounded by one. By
\lemref{dim:log-tails} and \propref{prop:unique}, there are constants
$K,c>0$ such that
\begin{equation}\label{eq:pressure-derivative-tail}
 Q_{\la,z}(x)\le e^{-ax/2}+2K e^{-c(x/4)^d}
 \qquad(\la\in[a,b],\ 0\le z\le1).
\end{equation}
The right-hand side is integrable against $x\,\nu_d(\mathrm dx)$.
By \propref{prop:unique},
$Q_{\la,z}(x)\to Q_{\la,0}(x)$ for each $\la,x$. Since
$p(\la,z)\to p(\la)$ by \propref{prop:perfect}, dominated convergence in
\eqref{eq:soft-pressure-increment} gives
\begin{equation}\label{eq:hard-pressure-increment}
 p(b)-p(a)=-\int_a^b u(\la)\,\mathrm d\la,\qquad
 u(\la)=\int_0^\infty x\,Q_{\la,0}(x)\,\nu_d(\mathrm dx).
\end{equation}

We next show that $u$ is continuous. Let $\la_j\to\la>0$ and let
$\mathsf d$ be a metric inducing weak convergence on probability laws on
$(0,\infty)$. By \propref{prop:unique}, we can choose $z_j\downarrow0$
such that
\[
 \mathsf d\bigl(\mathcal L(\cH_{\la_j,z_j}),
               \mathcal L(\cH_{\la_j,0})\bigr)<1/j.
\]
By \propref{prop:existence z = 0}, the laws of $\cH_{\la_j,z_j}$ are
relatively compact, and every subsequential limit solves the zero-activity
recursion at $\la$. Uniqueness implies convergence of these laws to
$\mathcal L(\cH_{\la,0})$. The metric bound gives the same convergence for
the laws of $\cH_{\la_j,0}$. Joint continuity of $f_\la(x,h,k)$ in
$(\la,h,k)$ and weak convergence of the independent message pairs imply
$Q_{\la_j,0}(x)\to Q_{\la,0}(x)$ for each $x$. Choosing $[a,b]$ to contain
$\la$ and all sufficiently large $\la_j$, we can apply
\eqref{eq:pressure-derivative-tail} and dominated convergence to conclude
that $u(\la_j)\to u(\la)$. The fundamental theorem of calculus applied to
\eqref{eq:hard-pressure-increment} now gives $p'(\la)=-u(\la)$ for every
$\la>0$.

Finally, fix $\la>0$. Let $(\xi_i)_{i\ge1}$ be a Poisson point process
with intensity $\nu_d$, and let $(H_i)_{i\ge1}$ be independent marks of
law $\mathcal L(\cH_{\la,0})$, independent of the point process. Define
\[
 a_i=\frac{e^{-\la\xi_i}H_i}{\sum_{j\ge1}e^{-\la\xi_j}H_j},
 \qquad
 A_i=\left(\sum_{j\ne i}e^{-\la\xi_j}H_j\right)^{-1}.
\]
The coefficients $a_i$ sum to one almost surely and satisfy
$a_i=f_\la(\xi_i,A_i,H_i)$. In the deletion form of Mecke's identity
\eqref{eq:mecke}, the remaining process supplies a message of law
$\mathcal L(\cH_{\la,0})$ by the recursion, and the inserted point carries
an independent mark of the same law. Consequently,
\[
 1=\E{\sum_{i\ge1}a_i}
 =\int_0^\infty Q_{\la,0}(x)\,\nu_d(\mathrm dx).
\]
For fixed $h,k$, we have
\[
 -\partial_x f_\la(x,h,k)
 =\la f_\la(x,h,k)\bigl(1-f_\la(x,h,k)\bigr).
\]
The derivative is bounded in absolute value by $\la/4$. Differentiation
under the expectation therefore gives
$-\partial_x Q_{\la,0}(x)\le\la Q_{\la,0}(x)$. The boundary term
$x^dQ_{\la,0}(x)$ vanishes at zero since $Q_{\la,0}\le1$, and at infinity
by \eqref{eq:pressure-derivative-tail}. Using
$\nu_d(\mathrm dx)=\mathrm d(x^d)$ and integrating by parts yields
\[
 1=-\int_0^\infty x^d\,\partial_x Q_{\la,0}(x)\,\mathrm dx
 \le\la\int_0^\infty x^d Q_{\la,0}(x)\,\mathrm dx
 =\frac{\la}{d}\,u(\la).
\]
This proves the lower bound.
\end{proof}

\section{Proof of Corollary \ref{cor:cost-budget}}\label{sec:corol}
\begin{proof}
Work on the probability-one event of pressure convergence at positive
rational inverse temperatures. Monotonicity of $p_n$ and continuity of
$p$ extend convergence to all $\la>0$. For
$m_n=\min_{M\in\mathscr M_n}\cC_n(M)$ and $\la>1$, retaining the
minimizer and comparing the partition functions at $\la$ and $1$ give
\[
 -p_n(\la)/\la \le m_n \le [p_n(1)-p_n(\la)]/(\la-1).
\]
By convexity and monotonicity the secant slopes of $p$ from $1$ have a
finite limit at infinity. Taking $n\to\infty$ and then $\la\to\infty$
in the bounds gives $m_n\to\alpha_d^\star$.

For the count, we use the exponential-tilting proof of the
G\"artner--Ellis lower bound~\cite{Gartner,Ellis}, as presented in
\cite[Theorem~2.3.6(b)]{DemboZeitouni}.
For $b>\alpha_d^\star$ write $N_n(b)=\#\{M\in\mathscr M_n:\cC_n(M)\le b\}$
and $s(b)=\inf_{t>0}\{p(t)+tb\}$. By \lemref{lem:pressure-derivative},
$u(\la)\to\infty$ as $\la\downarrow0$, and by convexity
$u(\la)\to\alpha_d^\star$ as $\la\to\infty$. Thus for any
$\alpha_d^\star<b'<b$ there is $\la>0$ with $u(\la)=b'$, and
$p(\la)+\la b'=s(b')$. With $0<\eta<b-b'$ and
$A_n=\{M:|\cC_n(M)-b'|\le\eta\}$, the Gibbs weights gives for all $t > 0$
\[
 e^{n[p_n(\la)+\la(b'-\eta)]}\mu_{\la,n}(A_n)
 \le N_n(b) \le e^{n[p_n(t)+tb]}.
\]
In the first inequality, we used the fact that $\cC_n(M) \geq b^\prime - \eta$ on $A_n$, so the weight is at most $e^{-n[p_n(\la) + \la(b^\prime - \eta)]}$. The proof of \corref{cor:energy} gives $\mu_{\la,n}(A_n)\to1$. Taking logarithmic limits gives the lower bound $s(b')$ as
$\eta\downarrow0$ and the upper bound $s(b)$ by taking the infimum over
$t>0$. The function $s$ is finite and concave on $(\alpha_d^\star,\infty)$,
so $b'\uparrow b$ gives equality.
\end{proof}

\section{Conclusion and Discussion}
In this work, we rigorously prove the exact limiting form of the free energy for the bipartite matching problem at every temperature. We discussed subsequent implications of the free energy formulas, namely concentration of the optimal matching cost and geometry of the solution space near the ground states. However, several additional questions remain unanswered. With the techniques outlined in this paper, it should be possible to study the universality of fluctuations for the bipartite assignment problem in every pseudodimension $d > 0$. In particular, we believe similar techniques could subsequently prove overlap concentration and Gaussian fluctuations between two independent replicas from the Gibbs measure around a deterministic limit in $\la$, as well as find other CLT-type fluctuations for the free energy and minimal assignment cost for the $\la > 0$ and $\la \to \infty$ regimes. One might hope that the same types of results emerge for different distributions of $C$, such as Gaussian or Laplace (see \cite{mordant2021maxima, Mordant}). A much more ambitious question is to consider other optimization problems beyond the matching problem, such as traveling salesman problem (TSP), or study the equivalent Gibbs measures  $\mu_{\la,n}$ introduced in this paper in the context of bipartite matchings for fixed ambient dimension \cite{caracciolo_scaling_2014, ambrosio_optimal_2019, ambrosio_quadratic_2021, huesmann_asymptotics_2024}. We leave these fascinating questions for future work.

\section*{Acknowledgments and AI Usage}

The author used ChatGPT-6.1 and ChatGPT-6.2 Sol in developing and drafting the manuscript, and in preparing the Lean formalization. The author subsequently rewrote, worked through, and verified every argument in the manuscript and takes full responsibility for the contents of the article. A machine-checked Lean~4 proof of \thmref{dim:main} is available in \href{https://github.com/shkroby/free-energy-assignment-lean}{https://github.com/shkroby/free-energy-assignment-lean}.

\appendix

\section{Appendix}\label{app:uniqueness}\label{app:simplex}

\begin{proof}[Proof of \lemref{lem:simplex}]
Write $\om(q)=(\log(q(1-q)))'=(1-2q)/(q(1-q))$. If a positive probability vector $q$ has three coordinates at least $\eta>0$, then for every absolutely summable $v$ with $\sum_i v_i=0$,
\begin{equation}\label{eq:simplex-curvature}
 \sum_i \om(q_i)v_i^2\ge\frac{\eta}{4}\|v\|_1^2.
\end{equation}
Take $q_1=\max_i q_i$ and put $g(s)=s\om(s)=(1-2s)/(1-s)$. If $q_1\le1/2$, then $\om(q_1)\ge0$ and $q_i\le(1-\eta)/2$ for $i>1$, so $g(q_i)\ge\eta$. If $q_1>1/2$, put $r=1-q_1$. For $i>1$ we have $q_i\le r-\eta$, and $g'\le-1$ gives $g(q_i)\ge g(r)+\eta$. Here $g(r)>0$ and $\om(q_1)=-g(r)/r$. By Cauchy--Schwarz, $\sum_{i>1}v_i^2/q_i\ge v_1^2/r$, which absorbs the contribution of $\om(q_1)$. In either case,
\[
 \sum_i \om(q_i)v_i^2\ge\eta\sum_{i>1}\frac{v_i^2}{q_i}
 \ge\eta\biggl(\sum_{i>1}|v_i|\biggr)^2
 \ge\frac{\eta}{4}\|v\|_1^2.
\]
For infinitely many coordinates the Cauchy--Schwarz bounds follow by finite truncation.

Now let $v=a-\wt a$ and $q(s)=\wt a+sv$ for $0\le s\le1$. The fundamental theorem of calculus gives
\[
 \cD(a,\wt a)=\int_0^1\sum_i\om(q_i(s))v_i^2\,\dd s.
\]
Only the coordinates exceeding $1/2$ at either endpoint can give negative integrands. There are at most two, and their integrals are finite since their segments stay in $(0,1)$. Tonelli's theorem applies to the remaining terms, proving the identity also for countably many coordinates. Choose any three coordinates and let $\eta>0$ be the smallest of their values at the two endpoints. They are at least $\eta$ along the segment, so~\eqref{eq:simplex-curvature} proves strict positivity whenever $v\ne0$.
\end{proof}

\printbibliography
\end{document}